\documentclass{article}
\usepackage{subfiles}
\usepackage{ifthen}

\newif\ifpreprint
\preprinttrue

\newcommand{\TheTitle}{A globally convergent fast iterative shrinkage-thresholding algorithm with a new momentum factor for single and multi-objective convex optimization}

\ifpreprint
\usepackage[nohyperref,abbrvbib,preprint]{jmlr2e}
\else
\usepackage[nohyperref,abbrvbib]{jmlr2e}
\fi

\jmlrheading{}{}{}{}{}{}{Hiroki Tanabe, Ellen H. Fukuda, and Nobuo Yamashita}
\ShortHeadings{A globally convergent FISTA}{Tanabe, Fukuda, and Yamashita}
\firstpageno{1}

\title{\TheTitle}
\author{\name Hiroki Tanabe \email tanabehiroki@amp.i.kyoto-u.ac.jp \\
       \name Ellen H. Fukuda \email ellen@i.kyoto-u.ac.jp \\
       \name Nobuo Yamashita \email nobuo@i.kyoto-u.ac.jp \\
       \addr Department of Applied Mathematics and Physics\\
       Graduate School of Informatics\\
       Kyoto University\\
       Yoshida-Honmachi, Sakyo-ku, Kyoto 606-8501, Japan
       }

\editor{}

\usepackage{url}
\usepackage{natbib}
\usepackage{graphicx}
\usepackage{adjustbox}
\usepackage{epstopdf}
\usepackage{setspace}
\usepackage{amsmath,amssymb}
\usepackage{mathtools}
\usepackage{cases}
\usepackage[section]{algorithm}
\usepackage[colorlinks=false,allbordercolors={1 1 1}]{hyperref}
\usepackage{algpseudocode}
\usepackage{derivative}
\usepackage{booktabs}
\usepackage{tabularray}
\usepackage{array}
\usepackage{here}
\usepackage{empheq}
\usepackage{siunitx}
\usepackage{csvsimple}
\usepackage{subcaption}
\usepackage{enumitem}
\usepackage[capitalize,nameinlink]{cleveref}
\let\etoolboxcsvloop\csvloop \let\csvloop\relax
\usepackage{autonum}
\let\csvloop\etoolboxcsvloop

\def\usebfsetcapital{\def\setcapital##1{\mathbf{##1}}}

\usebfsetcapital

\def\setR{\setcapital{R}}
\def\setRpos{\setcapital{R}_{\mathord{+}}}

\newcommand{\level}{\mathcal{L}}
\newcommand{\Beta}{\mathrm{B}}
\newcommand\condition[1]{\quad \text{#1}}
\newcommand\forallcondition[1]{\condition{for all~$#1$}}
\newcommand\eqand{\quad \text{and} \quad}

\DeclareMathOperator*{\argmin}{argmin}

\DeclareMathOperator{\interior}{int}

\DeclareMathOperator{\dom}{dom}

\DeclareMathOperator{\indicator}{\chi}
\DeclarePairedDelimiter{\abs}{\lvert}{\rvert}
\DeclarePairedDelimiter{\norm}{\lVert}{\rVert}

\DeclarePairedDelimiter{\set}{\lbrace}{\rbrace}
\DeclarePairedDelimiterX{\Set}[2]{\lbrace}{\rbrace}{#1\mathrel{}\delimsize\vert\mathrel{}#2}
\DeclarePairedDelimiterX{\innerp}[2]{\langle}{\rangle}{#1, #2}
\newcommand{\T}{\top\hspace{-1pt}}

\newcommand{\acc}{\mathrm{acc}}
\newcommand{\ra}[1]{\renewcommand{\arraystretch}{#1}}
\newtheorem{assumption}{Assumption}[section]

\crefname{equation}{}{}
\Crefname{equation}{Eq.}{Eqs.}
\crefname{enumi}{}{}
\crefname{figure}{Figure}{Figures}
\crefname{assumption}{Assumption}{Assumptions}
\crefname{line}{line}{lines}

\setlist[enumerate]{
    label=(\roman*)
}
\AtBeginEnvironment{theorem}{
    \setlist[enumerate]{
        label=(\roman*),
        ref=Theorem~\thetheorem~(\roman*)
    }
}
\AtEndEnvironment{theorem}{
    \setlist[enumerate]{
        label=(\roman*),
        ref=(\roman*)
    }
}
\AtBeginEnvironment{lemma}{
    \setlist[enumerate]{
        label=(\roman*),
        ref=Lemma~\thelemma~(\roman*)
    }
}
\AtEndEnvironment{lemma}{
    \setlist[enumerate]{
        label=(\roman*),
        ref=(\roman*)
    }
}
\AtBeginEnvironment{proposition}{
    \setlist[enumerate]{
        label=(\roman*),
        ref=Proposition~\theproposition~(\roman*)
    }
}
\AtEndEnvironment{proposition}{
    \setlist[enumerate]{
        label=(\roman*),
        ref=(\roman*)
    }
}
\AtBeginEnvironment{corollary}{
    \setlist[enumerate]{
        label=(\roman*),
        ref=Corollary~\thecorollary~(\roman*)
    }
}
\AtEndEnvironment{corollary}{
    \setlist[enumerate]{
        label=(\roman*),
        ref=(\roman*)
    }
}
\AtBeginEnvironment{remark}{
    \setlist[enumerate]{
        label=(\roman*),
        ref=Remark~\theremark~(\roman*)
    }
}
\AtEndEnvironment{remark}{
    \setlist[enumerate]{
        label=(\roman*),
        ref=(\roman*)
    }
}
\AtBeginEnvironment{definition}{
    \setlist[enumerate]{
        label=(\roman*),
        ref=Definition~\thedefinition~(\roman*)
    }
}
\AtEndEnvironment{definition}{
    \setlist[enumerate]{
        label=(\roman*),
        ref=(\roman*)
    }
}
\AtBeginEnvironment{assumption}{
    \setlist[enumerate]{
        label=(\roman*),
        ref=Assumption~\theassumption~(\roman*)
    }
}
\AtEndEnvironment{assumption}{
    \setlist[enumerate]{
        label=(\roman*),
        ref=(\roman*)
    }
}
\AtBeginEnvironment{example}{
    \setlist[enumerate]{
        label=(\roman*),
        ref=Example~\theexample~(\roman*)
    }
}
\AtEndEnvironment{example}{
    \setlist[enumerate]{
        label=(\roman*),
        ref=(\roman*)
    }
}
\newcounter{subcreftmpcnt}
\newcommand\romansubformat[1]{(\roman{#1})}
\makeatletter
\newcommand\subcref[2][\romansubformat]{
\ifcsname r@#2@cref\endcsname
  \cref@getcounter {#2}{\mylabel}
  \setcounter{subcreftmpcnt}{\mylabel}
  \hyperref[#2]{part~#1{subcreftmpcnt}}
 \else ?? \fi}
\newcommand\sublabelcref[2][\romansubformat]{
\ifcsname r@#2@cref\endcsname
  \cref@getcounter {#2}{\mylabel}
  \setcounter{subcreftmpcnt}{\mylabel}
  \hyperref[#2]{#1{subcreftmpcnt}}
 \else ?? \fi}
\newcommand\subCref[2][\romansubformat]{
\ifcsname r@#2@cref\endcsname
  \cref@getcounter {#2}{\mylabel}
  \setcounter{subcreftmpcnt}{\mylabel}
  \hyperref[#2]{Part~#1{subcreftmpcnt}}
 \else ?? \fi}

\makeatother

\begin{document}
\maketitle
\begin{abstract}
    Convex-composite optimization, which minimizes an objective function represented by the sum of a differentiable function and a convex one, is widely used in machine learning and signal/image processing.
    Fast Iterative Shrinkage Thresholding Algorithm (FISTA) is a typical method for solving this problem and has a global convergence rate of~$O(1 / k^2)$.
    Recently, this has been extended to multi-objective optimization, together with the proof of the~$O(1 / k^2)$ global convergence rate.
    However, its momentum factor is classical, and the convergence of its iterates has not been proven.
    In this work, introducing some additional hyperparameters~$(a, b)$, we propose another accelerated proximal gradient method with a general momentum factor, which is new even for the single-objective cases.
    We show that our proposed method also has a global convergence rate of~$O(1/k^2)$ for any $(a,b)$, and further that the generated sequence of iterates converges to a weak Pareto solution when~$a$ is positive, an essential property for the finite-time manifold identification.
    Moreover, we report numerical results with various~$(a,b)$, showing that some of these choices give better results than the classical momentum factors.
\end{abstract}

\ifpreprint
\else
\fi

\newboolean{isMain}
\setboolean{isMain}{true}

\section{Introduction} \label{sec: intro}
We consider the following convex-composite single~($m = 1$) or multi-objective~($m \ge 2$) optimization problem:
\[ \label{eq:MOP}
    \min_{x \in \setR^n} \quad F(x)
,\]
where~$F \colon \setR^n \to (\setR \cup \set*{\infty})^m$ is a vector-valued function with~$F \coloneqq (F_1, \dots, F_m)^\T$.
We assume that each component~$F_i \colon \setR^n \to \setR \cup \set*{\infty}$ is given by 
\[
    F_i(x) \coloneqq f_i(x) + g_i(x) \forallcondition{i = 1, \dots, m}
\]
with convex and continuously differentiable functions~$f_i \colon \setR^n \to \setR, i = 1, \dots, m$ and closed, proper and convex functions~$g_i \colon \setR^n \to \setR \cup \set*{\infty}, i = 1, \dots, m$, and each~$\nabla f_i$ is Lipschitz continuous.
As suggested in~\cite{Tanabe2019}, this problem involves many important classes.
For example, it can express a convex-constrained problem if each~$g_i$ is the indicator function of a convex set~$S$, i.e.,
\[ \label{eq:indicator}
    \indicator_S(x) \coloneqq \begin{dcases}
        0, & \text{if } x \in S,\\
        \infty, & \text{otherwise}
    .\end{dcases}
\] 

Multi-objective optimization has many applications in engineering~\citep{Eschenauer1990}, statistics~\citep{Carrizosa1998}, and machine learning (particularly multi-task learning~\citep{Sener2018,Lin2019} and neural architecture search~\citep{Kim2017,Dong2018,Elsken2019}).
In the multi-objective case, no single point minimizes all objective functions simultaneously in general.
Therefore, we use the concept of \emph{Pareto optimality}.
We call a point weakly Pareto optimal if there is no other point where the objective function values are strictly smaller.
This generalizes the usual optimality for single-objective problems.
In other words, single-objective problems are considered to be included in multi-objective ones.
Hence, in the following, unless otherwise noted, we refer to~\cref{eq:MOP} as multi-objective, including the case where~$m = 1$.

One of the main strategies for multi-objective problems is the~\emph{scalarization approach}~\citep{Gass1955,Geoffrion1968,Zadeh1963}, which reduces the original multi-objective problem into a parameterized (or weighted) scalar-valued problem.
However, it requires an \emph{a priori} parameters (or weights) selection, which might be challenging.
The meta-heuristics~\citep{Gandibleux2004} is also popular, but it has no theoretical convergence properties under reasonable assumptions.

Many descent methods have been developed in recent years~\citep{Fukuda2014}, overcoming those drawbacks.
They decrease all objective values simultaneously at each iteration, and their global convergence property can be analyzed under reasonable assumptions.
For example, the steepest descent method~\citep{Fliege2000,Fliege2019,Desideri2012} converges globally to Pareto solutions for differentiable multi-objective problems.
From a practical point of view, its applicability has also been reported in multi-task learning~\citep{Sener2018,Lin2019}.
Afterwards, the projected gradient~\citep{Fukuda2013}, Newton's~\citep{Fliege2009,Goncalves2021}, trust-region~\citep{Carrizo2016}, and conjugate gradient methods~\citep{LucambioPerez2018} were also considered.
Moreover, the proximal point~\citep{Bonnel2005} and the inertial forward-backward methods~\citep{Bot2018} can solve infinite-dimensional vector optimization problems.

For~\cref{eq:MOP}, the proximal gradient method~\citep{Tanabe2019,Tanabe2022} is effective.
Using it, the merit function~\citep{Tanabe2022b}, which returns zero at the Pareto solutions and strictly positive values otherwise, converges to zero with rate~$O(1 / k)$ under reasonable assumptions.
It is also shown that the generated sequence of iterates converges to a weak Pareto solution~\citep{Bello-Cruz2022}.
On the other hand, the accelerated proximal gradient method~\citep{Tanabe2022a}, which generalizes the Fast Iterative Shrinkage Thresholding Algorithm (FISTA)~\citep{Beck2009} for convex-composite single-objective problems, has also been considered, along with a proof of the merit function's~$O(1 / k^2)$ convergence rate.
However, the momentum factor used there is classical ($t_1 = 1, t_{k + 1} = \sqrt{t_k^2 + 1 / 4} + 1 / 2$), and the iterates' convergence is not proven.

This paper generalizes the associated factor by~$t_1 = 1, t_{k + 1} = \sqrt{t_k^2 - a t_k + b} + 1 / 2$ with hyperparameters~$a \in [0, 1), b \in [a^2 / 4, 1 / 4]$.
This is new even in the single-objective context, and it generalizes well-known factors.
For example, when~$a = 0$ and~$b = 1 / 4$, it reduces to~$t_1 = 1, t_{k + 1} = \sqrt{t_k^2 + 1 / 4} + 1 / 2$, proposed in~\cite{Nesterov1983,Beck2009}, and when~$b = a^2 / 4$, it gives~$t_k = (1 - a) k / 2 + (1 + a) / 2$, suggested in~\cite{Chambolle2015,Attouch2016,Attouch2018,Su2016}.
We show that the merit function converges to zero with rate~$O(1 / k^2)$ for any~$(a, b)$.
In addition, we prove the iterates' convergence to a weak Pareto solution when~$a > 0$.
As discussed in \cref{sec: convergence sequence}, this suggests that the proposed method might achieve finite-iteration manifold (active set) identification~\citep{Sun2019} without the assumption of strong convexity.
Furthermore, we carry out numerical experiments with various~$(a, b)$ and observe that some~$(a, b)$ yield better results than the classical factors.

The outline of this paper is as follows.
We present some notations and definitions used in this paper in \cref{sec:notations}.
\Cref{sec: acc prox} recalls the accelerated proximal gradient method for~\cref{eq:MOP} and its associated results.
We generalize the momentum factor and prove that it preserves an~$O(1 / k^2)$ convergence rate in \cref{sec: generalization}, and we demonstrate the convergence of the iterates in \cref{sec: convergence sequence}.
Finally, \cref{sec: experiments} provides numerical experiments and compares the numerical performances depending on the hyperparameters.

\ifthenelse{\boolean{isMain}}{ 
}{
   {\bibliographystyle{jorsj} 
   \bibliography{library}} 
}

\section{Preliminaries} \label{sec: preliminaries}
\subsection{Definitions and notations} \label{sec:notations}
For every natural number~$d$, write the~$d$-dimensional real space by~$\setR^d$, and define
\[
    \setRpos^d \coloneqq \Set*{v \in \setR^d}{v_i \ge 0, i = 1, \dots, d}
.\]
This induces the partial orders: for any~$v^1, v^2 \in \setR^d$,~$v^1 \le v^2$ (alternatively, $v^2 \ge v^1$) if~$v^2 - v^1 \in \setRpos^d$ and~$v^1 < v^2$ (alternatively,~$v^2 > v^1$) if~$v^2 - v^1 \in \interior \setRpos^d$.
In other words, $v^1 \le v^2$ and~$v^1 < v^2$ mean that~$v^1_i \le v^2_i$ and~$v^1_i < v^2_i$ for all~$i = 1, \dots, d$, respectively.
Furthermore, let~$\innerp*{\cdot}{\cdot}$ be the Euclidean inner product in~$\setR^d$, i.e.,~$\innerp*{v^1}{v^2} \coloneqq \sum_{i = 1}^d {v^1_i v^2_i}$, and let~$\norm*{ \cdot }$ be the Euclidean norm, i.e., $\norm*{v} \coloneqq \sqrt{\innerp*{v}{v}}$.
Moreover, we define the $\ell_1$-norm and~$\ell_\infty$-norm by~$\norm*{v}_1 \coloneqq \sum_{i = 1}^{m} \abs*{v_i}$ and~$\norm*{v}_\infty \coloneqq \max_{i = 1, \dots, d} \abs*{v_i}$, respectively.

We introduce some concepts used in the problem~\cref{eq:MOP}.
Recall that
\[ \label{eq:weak Pareto}
    X^\ast \coloneqq \Set*{x^\ast \in \setR^n}{\text{There does not exist~$x \in \setR^n$ such that~$F(x) < F(x^\ast)$}}
\]
is the set of \emph{weakly Pareto optimal} solutions for~\cref{eq:MOP}.
When~$m = 1$,~$X^\ast$ reduces to the optimal solution set.
Moreover, define the effective domain of~$F$ by
\[
    \dom F \coloneqq \Set*{x \in \setR^n}{F(x) < \infty}
,\]
and write the level set of~$F$ on~$c \in \setR^m$ as
\[ \label{eq:level set}
    \level_F(c) \coloneqq \Set*{x \in \setR^n}{F(x) \le c}
.\]
Furthermore, we express the image of~$A \subseteq \setR^n$ and the inverse image of~$B \subseteq (\setR \cup \set*{\infty})^m$ under~$F$ as
\[
    F(A) \coloneqq \Set*{F(x) \in \setR^m}{x \in A} \eqand F^{-1}(B) \coloneqq \Set*{x \in \setR^n}{F(x) \in B},
\]
respectively.

Finally, let us recall the merit function~$u_0 \colon \setR^n \to \setR \cup \set*{\infty}$ for~\cref{eq:MOP} proposed in~\cite{Tanabe2022b}:
\[ \label{eq:u_0}
    u_0(x) \coloneqq \sup_{z \in \setR^n} \min_{i = 1, \dots, m} [ F_i(x) - F_i(z) ]
,\]
which returns zero at optimal solutions and strictly positive values otherwise.
The following theorem shows that $u_0$ is a merit function in the Pareto sense.
\begin{theorem}~\cite[Theorem 3.1]{Tanabe2022b} \label{thm:merit Pareto}
    Let~$u_0$ be defined by~\cref{eq:u_0}.
    Then,~$u_0(x) \ge 0$ for all~$x \in \setR^n$.
    Moreover,~$x \in \setR^n$ is weakly Pareto optimal for~\cref{eq:MOP} if and only if~$u_0(x) = 0$.
\end{theorem}
Note that when~$m = 1$, we have
\[
    u_0(x) = F_1(x) - F_1^\ast
,\]
where~$F_1^\ast$ is the optimal objective value.
Clearly, this is a merit function for scalar-valued optimization.

\subsection{The accelerated proximal gradient method for multi-objective optimization} \label{sec: acc prox}
This subsection recalls the accelerated proximal gradient method for~\cref{eq:MOP} proposed in~\cite{Tanabe2022a} and its main results.
Recall that each~$F_i$ is the sum of a continuously differentiable function~$f_i$ and a closed, proper, and convex function~$g_i$, and that~$\nabla f_i$ is Lipschitz continuous with Lipschitz constant~$L_i > 0$.
Define
\[ \label{eq:L}
    L \coloneqq \max_{i = 1, \dots, m} L_i
.\] 
The method solves the following subproblem at each iteration for given~$x \in \dom F$,~$y \in \setR^n$, and~$\ell \ge L$:
\[ \label{eq:acc prox subprob}
    \min_{z \in \setR^n} \quad \varphi^\acc_\ell(z; x, y) 
,\]
where
\[ \label{eq:varphi acc}
\varphi^\acc_\ell(z; x, y) \coloneqq \max_{i = 1, \dots, m} \left[ \innerp*{\nabla f_i(y)}{z - y} + g_i(z) + f_i(y) - F_i(x) \right] + \frac{\ell}{2} \norm*{ z - y }^2
.\]
From the strong convexity,~\cref{eq:acc prox subprob} has a unique optimal solution~$p^\acc_\ell(x, y)$, i.e.,
\[ \label{eq:p theta acc}
    p^\acc_\ell(x, y) \coloneqq \argmin_{z \in \setR^n} \varphi^\acc_\ell(z; x, y)
.\]
The following proposition characterizes weak Pareto optimality in terms of the mapping~$p^\acc_\ell$.
\begin{proposition}~\cite[Proposition 4.1 (i)]{Tanabe2022a} \label{thm:acc prox termination}
    Let~$p^\acc_\ell(x, y)$ be defined by~\cref{eq:p theta acc}.
    Then,~$y \in \setR^n$ is weakly Pareto optimal for~\cref{eq:MOP} if and only if~$p^\acc_\ell(x, y) = y$ for some~$x \in \setR^n$.
\end{proposition}
This implies that using~$\norm*{p^\acc_\ell(x, y) - y}_\infty < \varepsilon$ for some~$\varepsilon > 0$ is reasonable as the stopping criteria.
We state below the accelerated proximal gradient method for~\cref{eq:MOP}.
\begin{algorithm}[hbtp]
    \caption{Accelerated proximal gradient method for~\cref{eq:MOP}}
    \label{alg:acc-pgm}
    \begin{algorithmic}[1]
        \Require Set~$x^0 = y^1 \in \dom F, \ell \ge L, \varepsilon > 0$.
        \Ensure $x^\ast$: A weakly Pareto optimal point
        \State $k \gets 1$
        \State $t_1 \gets 1$ \label{line:t ini}
        \While{$\norm*{p^\acc_\ell(x^{k - 1}, y^k) - y^k}_\infty \ge \varepsilon$}
        \State $x^k \gets p^\acc_\ell(x^{k - 1}, y^k)$
        \State $t_{k + 1} \gets \sqrt{t_k^2 + 1 / 4} + 1/2$ \label{line:t rr}
        \State $\gamma_k \gets (t_k - 1) / t_{k + 1}$ \label{line:gamma}
        \State $y^{k + 1} \gets x^k + \gamma_k (x^k - x^{k - 1})$ \label{line:y}
        \State $k \gets k + 1$
        \EndWhile
    \end{algorithmic}
\end{algorithm}

\Cref{alg:acc-pgm} generates~$\set*{x^k}$ such that~$\set*{u_0(x^k)}$ converges to zero with rate~$O(1 / k^2)$ under the following assumption.
This assumption is also used to analyze the proximal gradient method without acceleration~\citep{Tanabe2022} and is not particularly strong as suggested in~\citep[Remark 5.2]{Tanabe2022}.
\begin{assumption}~\cite[Assumption 5.1]{Tanabe2022} \label{asm:bound}
    Let~$X^\ast$ and~$\level_F$ be defined by~\cref{eq:weak Pareto,eq:level set}, respectively.
    Then, for all~$x \in \level_F(F(x^0))$, there exists~$x^\ast \in X^\ast$ such that~$F(x^\ast) \le F(x)$ and
    \[ \label{eq:R}
        R \coloneqq \sup_{F^\ast \in F(X^\ast \cap \level_F(F(x^0))} \inf_{z \in F^{-1}(\set*{F^\ast})} \norm*{z - x^0}^2 < \infty
    .\] 
\end{assumption}

\begin{theorem}~\cite[Theorem 5.2]{Tanabe2022a} \label{thm:conv rate}
    Under \cref{asm:bound}, \cref{alg:acc-pgm} generates~$\set*{x^k}$ such that
    \[
        u_0(x^k) \le \frac{2 \ell R}{(k + 1)^2} \forallcondition{k \ge 1}
    ,\] 
    where~$R \ge 0$ is given by~\cref{eq:R}, and~$u_0$ is a merit function defined by~\cref{eq:u_0}.
\end{theorem}

The following corollary shows the global convergence of \cref{alg:acc-pgm}.
\begin{corollary}~\cite[Corollary 5.2]{Tanabe2022a} \label{thm:accumulation point}
    Suppose that \cref{asm:bound} holds.
    Then, every accumulation point of~$\set*{x^k}$ generated by \cref{alg:acc-pgm} is weakly Pareto optimal for~\cref{eq:MOP}.
\end{corollary}

\ifthenelse{\boolean{isMain}}{ 
}{
   {\bibliographystyle{jorsj} 
   \bibliography{library}} 
}

\section{Generalization of the momentum factor and convergence rate analysis} \label{sec: generalization}
This section generalizes the momentum factor~$\set*{t_k}$ used in \cref{alg:acc-pgm} and shows that the~$O(1 / k^2)$ convergence rate also holds in that case.
First, we describe below the algorithm in which we replace \cref{line:t rr} of \cref{alg:acc-pgm} by a formula using given constants~$a \in [0, 1)$ and~$b \in [a^2 / 4, 1 / 4]$:
\begin{algorithm}[hbtp]
    \caption{Accelerated proximal gradient method with general stepsizes for~\cref{eq:MOP}}
    \label{alg:acc-pgm general}
    \begin{algorithmic}[1]
        \Require Set~$x^0 = y^1 \in \dom F, \ell \ge L, \varepsilon > 0, a \in [0, 1), b \in [a^2 / 4, 1 / 4]$.
        \Ensure $x^\ast$: A weakly Pareto optimal point
        \State $k \gets 1$
        \State $t_1 \gets 1$ \label{line:t ini general}
        \While{$\norm*{p^\acc_\ell(x^{k - 1}, y^k) - y^k}_\infty \ge \varepsilon$}
        \State $x^k \gets p^\acc_\ell(x^{k - 1}, y^k)$
        \State $t_{k + 1} \gets \sqrt{t_k^2 - a t_k + b} + 1/2$ \label{line:t rr general}
        \State $\gamma_k \gets (t_k - 1) / t_{k + 1}$ \label{line:gamma general}
        \State $y^{k + 1} \gets x^k + \gamma_k (x^k - x^{k - 1})$ \label{line:y general}
        \State $k \gets k + 1$
        \EndWhile
    \end{algorithmic}
\end{algorithm}

The sequence~$\set*{t_k}$ defined in~\cref{line:t ini general,line:t rr general} of \cref{alg:acc-pgm} generalizes the well-known momentum factors in single-objective accelerated methods.
For example, when~$a = 0$ and~$b = 1 / 4$, they coincide with the one in \cref{alg:acc-pgm} and the original FISTA~\citep{Nesterov1983,Beck2009} ($t_1 = 1$ and~$t_{k + 1} = (1 + \sqrt{1 + 4 t_k^2}) / 2$).
Moreover, if~$b = a^2 / 4$, then~$\set*{t_k}$ has the general term~$t_k = (1 - a) k / 2 + (1 + a) / 2$, which corresponds to the one used in~\cite{Chambolle2015,Su2016,Attouch2016,Attouch2018}.
This means that our generalization allows a finer tuning of the algorithm by varying~$a$ and~$b$.

We present below the main theorem of this section.
\begin{theorem} \label{thm:main theorem of stepsize section}
    Let~$\set*{x^k}$ be a sequence generated by \cref{alg:acc-pgm general} and recall that~$u_0$ is given by~\cref{eq:u_0}.
    Then, the following two equations hold:
    \begin{enumerate}
    \item $F_i(x^k) \le F_i(x^0)$ for all~$i = 1, \dots, m$ and~$k \ge 0$; \label{thm:main theorem of stepsize section:less than x0}
    \item $u_0(x^k) = O(1 / k^2)$ as $k \to \infty$ under \cref{asm:bound}. \label{thm:main theomrem of stepsize section:O(1/k2)}
    \end{enumerate}
\end{theorem}
\subCref{thm:main theorem of stepsize section:less than x0} means that~$\set*{x^k} \subseteq \level_F(F(x^0))$, where~$\level_F$ denotes the level set of~$F$ (cf.~\cref{eq:level set}).
Note, however, that the objective functions are generally not monotonically non-increasing.
\subCref{thm:main theomrem of stepsize section:O(1/k2)} also claims the global convergence rate.

Before proving \cref{thm:main theorem of stepsize section}, let us give several lemmas.
First, we present some properties of~$\set*{t_k}$ and~$\set*{\gamma_k}$.
\begin{lemma} \label{thm:t}
    Let~$\set*{t_k}$ and~$\set*{\gamma_k}$ be defined by \cref{line:t ini,line:t rr,line:gamma} in \cref{alg:acc-pgm} for arbitrary~$a \in [0, 1)$ and~$b \in [a^2 / 4, 1 / 4]$.
    Then, the following inequalities hold for all~$k \ge 1$.
    \begin{enumerate}
        \item $t_{k + 1} \ge t_k + \dfrac{1 - a}{2}$ and~$t_k \ge \dfrac{1 - a}{2} k + \dfrac{1 + a}{2}$; \label{thm:t:t geq}
        \item $t_{k + 1} \le t_k + \dfrac{1 - a + \sqrt{4b - a^2}}{2}$ and~$t_k \le \dfrac{1 - a + \sqrt{4b - a^2}}{2} (k - 1) + 1 \le k$; \label{thm:t:t leq}
        \item $t_k^2 - t_{k + 1}^2 + t_{k + 1} = a t_k - b + \dfrac{1}{4} \ge a t_k$; \label{thm:t:t over-relax geq}
        \item $0 \le \gamma_k \le \dfrac{k - 1}{k + 1 / 2}$; \label{thm:t:gamma}
        \item $1 - \gamma_k^2 \ge \dfrac{1}{t_k}$. \label{thm:t:t moment}
    \end{enumerate}
\end{lemma}
\begin{proof}
    \sublabelcref{thm:t:t geq}:
    From the definition of~$\set*{t_k}$, we have
    \[ \label{eq:t cs} 
        t_{k + 1} = \sqrt{t_k^2 - a t_k + b} + \frac{1}{2}
        = \sqrt{\left( t_k - \frac{a}{2} \right)^2 + \left( b - \frac{a^2}{4} \right)} + \frac{1}{2}
    .\]
    Since~$b \ge a^2 / 4$, we get
    \[
        t_{k + 1} \ge \abs*{t_k - \frac{a}{2}} + \frac{1}{2}
    .\]
    Since~$t_1 = 1 \ge a / 2$, we can quickly see that~$t_k \ge a / 2$ for any~$k$ by induction.
    Thus, we have
    \[
        t_{k + 1} \ge t_k + \frac{1 - a}{2}
    .\]
    Applying the above inequality recursively, we obtain
    \[
        t_k \ge \frac{1 - a}{2} (k - 1) + t_1 = \frac{1 - a}{2} k + \frac{1 + a}{2}
    .\]

    \sublabelcref{thm:t:t leq}:
    From~\cref{eq:t cs} and the relation~$\sqrt{\alpha + \beta} \le \sqrt{\alpha} + \sqrt{\beta}$ with~$\alpha, \beta \ge 0$, we get the first inequality.
    Using it recursively, it follows that
    \[
        t_k \le \frac{1 - a + \sqrt{4 b - a^2}}{2} (k - 1) + t_1 = \frac{1 - a + \sqrt{4 b - a^2}}{2} (k - 1) + 1
    .\] 
    Since~$a \in [0, 1), b \in [a^2 / 4, 1 / 4]$, we observe that
    \[
        \frac{1 - a + \sqrt{4 b - a^2}}{2} \le \frac{1 - a + \sqrt{1 - a^2}}{2} \le 1
    .\] 
    Hence, the above two inequalities lead to the desired result.

    \sublabelcref{thm:t:t over-relax geq}:
    An easy computation shows that
    \[
        \begin{split}
            t_k^2 - t_{k + 1}^2 + t_{k + 1} &= t_k^2 - \left[ \sqrt{t_k^2 - a t_k + b} + \frac{1}{2} \right]^2 + \sqrt{t_k^2 - a t_k + b} + \frac{1}{2} \\
                                            &= a t_k - b + \frac{1}{4} \ge a t_k
    ,\end{split}
    \] 
    where the inequality holds since~$b \le 1 / 4$.

    \sublabelcref{thm:t:gamma}:
    The first inequlity is clear from the definition of~$\gamma_k$ since \subcref{thm:t:t geq} yields~$t_k \ge 1$.
    Again, the definition of~$\gamma_k$ and \subcref{thm:t:t geq} give
    \[
        \gamma_k = \frac{t_k - 1}{t_{k + 1}} \le \frac{t_k - 1}{t_k + (1 - a) / 2} = 1 - \frac{3 - a}{2 t_k + 1 - a}
    .\] 
    Combining with \subcref{thm:t:t leq}, we get
    \[ \label{eq:gamma}
        \begin{split}
            \gamma_k &\le 1 - \frac{3 - a}{\left(1 - a + \sqrt{4 b - a^2} \right)(k - 1) + 3 - a} \\
        &= \frac{\left( 1 - a + \sqrt{4 b - a^2} \right) (k - 1)}{\left(1 - a + \sqrt{4 b - a^2} \right)(k - 1) + 3 - a} \\
        &= \frac{k - 1}{k - 1 + (3 - a) / \left( 1 - a + \sqrt{4 b - a^2} \right)}
        .\end{split}
    \] 
    On the other hand, it follows that
    \[ \label{eq:min a b}
        \min_{a \in [0, 1), b \in [a^2 / 4, 1 / 4]} \frac{3 - a}{1 - a + \sqrt{4 b - a^2}} = \min_{a \in [0, 1)} \frac{3 - a}{1 - a + \sqrt{1 - a^2}} = \frac{3}{2}
    ,\]
    where the second equality follows from the monotonic non-decreasing property implied by
    \[
        \odv*{\left(\frac{3 - a}{1 - a + \sqrt{1 - a^2}}\right)}{a} = \frac{2 \sqrt{1 - a^2} + 3a - 1}{\left( \sqrt{1 - a^2} - a + 1 \right)^2 \sqrt{1 - a^2}} > 0 \forallcondition{a \in [0, 1)}
    .\]  
    Combining~\cref{eq:gamma,eq:min a b}, we obtain~$\gamma_k \le (k - 1) / (k + 1 / 2)$.

    \sublabelcref{thm:t:t moment}:
    \subCref{thm:t:t geq} implies that~$t_{k + 1} > t_k \ge 1$.
    Thus, the definition of~$\gamma_k$ implies that
    \[
        1 - \gamma_k^2 = 1 - \left( \frac{t_k - 1}{t_{k + 1}} \right)^2 \ge 1 - \left( \frac{t_k - 1}{t_k} \right)^2
        = \frac{2 t_k - 1}{t_k^2} \ge \frac{2 t_k - t_k}{t_k^2} = \frac{1}{t_k}
    .\]
\end{proof}

As in~\cite{Tanabe2022a}, we also introduce~$\sigma_k \colon \setR^n \to \setR \cup \set*{- \infty}$ and~$\rho_k \colon \setR^n \to \setR$ for~$k \ge 0$ as follows, which assist the analysis:
\[ \label{eq:sigma rho}
\begin{gathered} 
    \sigma_k(z) \coloneqq \min_{i = 1, \dots, m}\left[ F_i(x^k) - F_i(z) \right], \\
        \rho_k(z) \coloneqq \norm*{t_{k + 1} x^{k + 1} - (t_{k + 1} - 1) x^k - z}^2
.\end{gathered}
\]
The following lemma on~$\sigma_k$ is helpful in the subsequent discussions.
\begin{lemma}~\cite[Lemma 5.1]{Tanabe2022a} \label{thm:sigma}
    Let~$\set*{x^k}$ and~$\set*{y^k}$ be sequences generated by \cref{alg:acc-pgm general}.
    Then, the following inequalities hold for all~$z \in \setR^n$ and~$k \ge 0$:
    \begin{enumerate}
        \item \label{thm:sigma:1} $\begin{multlined}[t] \sigma_{k + 1}(z) \le - \frac{\ell}{2} \left( 2 \innerp*{x^{k + 1} - y^{k + 1}}{y^{k + 1} - z} + \norm*{x^{k + 1} - y^{k + 1}}^2 \right)\\
            \displaystyle - \frac{\ell - L}{2} \norm*{x^{k + 1} - y^{k + 1}}^2;\end{multlined}$
        \item \label{thm:sigma:2}$\begin{multlined}[t]
                \sigma_k(z) - \sigma_{k + 1}(z) \ge \frac{\ell}{2} \left( 2 \innerp*{x^{k + 1} - y^{k + 1}}{y^{k + 1} - x^k} + \norm*{x^{k + 1} - y^{k + 1}}^2 \right) \\
            + \frac{\ell - L}{2} \norm*{x^{k + 1} - y^{k + 1}}^2
        .\end{multlined}$
    \end{enumerate}
\end{lemma}

Therefore, from \cref{thm:t:t moment}, we can obtain the following result quickly in the same way as in the proof of~\cite[Corollary 5.1]{Tanabe2022a}.
\begin{lemma} \label{thm:sigma k1 k2}
    Let~$\set*{x^k}$ and~$\set*{y^k}$ be sequences generated by \cref{alg:acc-pgm general}.
    Then, we have
    \begin{multline}
        \sigma_{k_1}(z) - \sigma_{k_2}(z) \\
        \ge \frac{\ell}{2} \left( \norm*{x^{k_2} - x^{k_2 - 1}}^2 - \norm*{x^{k_1} - x^{k_1 - 1}}^2 + \sum_{k = k_1}^{k_2 - 1} \frac{1}{t_k} \norm*{x^k - x^{k - 1}}^2 \right)
    \end{multline}
    for any~$k_2 \ge k_1 \ge 1$.
\end{lemma}

We can now show the first part of \cref{thm:main theorem of stepsize section}.
\begin{proof}[Proof of \cref{thm:main theorem of stepsize section:less than x0}]
    From \cref{thm:sigma k1 k2}, we can prove this part with similar arguments used in the proof of~\cite[Theorem~5.1]{Tanabe2022a}.
\end{proof}

The next step is to prepare the proof of \cref{thm:main theomrem of stepsize section:O(1/k2)}.
First, we mention the following relation, used frequently hereafter:
\begin{align}
        &\norm*{v^2 - v^1}^2 + 2 \innerp*{v^2 - v^1}{v^1 - v^3} = \norm*{ v^2 - v^3 }^2 - \norm*{v^1 - v^3}^2 , \label{eq:Pythagoras} \\
        &\sum_{s=1}^r \sum_{p=1}^s A_p = \sum_{p=1}^r \sum_{s=p}^r A_p \label{eq:change sum}
\end{align}
for any vectors~$v^1, v^2, v^3$ and sequence~$\set*{A_p}$.
With these, we show the lemma below, which is similar to~\cite[Lemma 5.2]{Tanabe2022a} but more complex due to the generalization of~$\set*{t_k}$.
\begin{lemma} \label{thm:key relation}
    Let~$\set*{x^k}$ and~$\set*{y^k}$ be sequences generated by \cref{alg:acc-pgm general}.
    Also, let~$\sigma_k$ and~$\rho_k$ be defined by~\cref{eq:sigma rho}.
    Then, we have
    \begin{align}
        \MoveEqLeft \frac{\ell}{2} \norm*{x^0 - z}^2 \\
        \ge{}& \frac{1}{1 - a} \left[ t_{k + 1}^2 - a t_{k + 1} + \left( \frac{1}{4} - b \right) k \right] \sigma_{k + 1}(z) \\
        &+ \frac{\ell}{2 (1 - a)} \left[ a (t_{k + 1}^2 - t_{k + 1}) + \left( \frac{1}{4} - b \right) k \right] \norm*{x^{k + 1} - x^k}^2 \\
        &+ \frac{\ell}{2 (1 - a)} \sum_{p = 1}^{k} \left[ a^2 (t_p - 1) + \left( \frac{1}{4} - b \right) \frac{p - t_p + a(t_p - 1)}{t_p} \right] \norm*{x^p - x^{p - 1}}^2 \\
        &+ \frac{\ell}{2} \rho_k(z) + \frac{\ell - L}{2} \sum_{p = 1}^{k} t_{p + 1}^2 \norm*{x^{p + 1} - y^{p + 1}}^2
    \end{align}
    for all~$k \ge 0$ and~$z \in \setR^n$.
\end{lemma}
\begin{proof}
    Let~$p \ge 1$ and~$z \in \setR^n$.
    Recall that \cref{thm:sigma} gives
    \begin{gather}
        \begin{multlined}
        - \sigma_{p + 1}(z) \ge \frac{\ell}{2} \left[ 2 \innerp*{x^{p + 1} - y^{p + 1}}{y^{p + 1} - z} + \norm*{x^{p + 1} - y^{p + 1}}^2 \right] \\
        + \frac{\ell - L}{2} \norm*{x^{p + 1} - y^{p + 1}}^2,
        \end{multlined} \\
        \begin{multlined}
            \sigma_p(z) - \sigma_{p + 1}(z) \ge \frac{\ell}{2} \left[ 2 \innerp*{x^{p + 1} - y^{p + 1}}{y^{p + 1} - x^p} + \norm*{x^{p + 1} - y^{p + 1}}^2 \right] \\
            + \frac{\ell - L}{2} \norm*{x^{p + 1} - y^{p + 1}}^2
        .\end{multlined}
        \end{gather}
        We then multiply the second inequality above by $(t_{p + 1} - 1)$ and add it to the first one:
    \begin{multline}
        (t_{p + 1} - 1) \sigma_p(z) - t_{p + 1} \sigma_{p + 1}(z) \\
        \ge \frac{\ell}{2} \left[ t_{p + 1} \norm*{x^{p + 1} - y^{p + 1}}^2 + 2 \innerp*{x^{p + 1} - y^{p + 1}}{t_{p + 1} y^{p + 1} - (t_{p + 1} - 1)x^p - z} \right] \\
        + \frac{\ell - L}{2} t_{p + 1} \norm*{x^{p + 1} - y^{p + 1}}^2
    .\end{multline}
    Multiplying this inequality by~$t_{p + 1}$ and using the relation~$t_p^2 = t_{p + 1}^2 - t_{p + 1} + (a t_p - b + 1/4)$ (cf. \cref{thm:t:t over-relax geq}), we get
    \begin{multline}
        t_p^2 \sigma_p(z) - t_{p + 1}^2 \sigma_{p + 1}(z) \ge \frac{\ell}{2} \Bigl[ \norm*{t_{p + 1} (x^{p + 1} - y^{p + 1})}^2  \\
        + 2 t_{p + 1} \innerp*{x^{p + 1} - y^{p + 1}}{t_{p + 1} y^{p + 1} - (t_{p + 1} - 1)x^p - z} \Bigr] \\
        + \frac{\ell - L}{2} t_{p + 1}^2 \norm*{x^{p + 1} - y^{p + 1}}^2 + \left( a t_p - b + \frac{1}{4} \right) \sigma_p(z)
    .\end{multline}
    Applying~\cref{eq:Pythagoras} to the right-hand side of the last inequality with
    \[
        v^1 \coloneqq t_{p + 1} y^{p + 1}, \quad v^2 \coloneqq t_{p + 1} x^{p + 1}, \quad v^3 \coloneqq (t_{p + 1} - 1) x^p + z
    .\]
    we get
    \begin{multline}
        t_p^2 \sigma_p(z) - t_{p + 1}^2 \sigma_{p + 1}(z) \\
        \ge \frac{\ell}{2} \left[ \norm*{t_{p + 1} x^{p + 1} - (t_{p + 1} - 1) x^p - z}^2 - \norm*{t_{p + 1} y^{p + 1} - (t_{p + 1} - 1) x^p - z}^2 \right] \\
        + \frac{\ell - L}{2} t_{p + 1}^2 \norm*{x^{p + 1} - y^{p + 1}}^2 + \left( a t_p - b + \frac{1}{4} \right) \sigma_p(z)
    .\end{multline}
    Recall that~$\rho_p(z) \coloneqq \norm*{t_{p + 1} x^{p + 1} - (t_{p + 1} - 1) x^p - z}^2$.
    Then, considering the definition of~$y^p$ given in \cref{line:y} of \cref{alg:acc-pgm}, we obtain
    \begin{multline}
        t_p^2 \sigma_p(z) - t_{p + 1}^2 \sigma_{p + 1}(z) \\
        \ge \frac{\ell}{2} \left[ \rho_p(z) - \rho_{p - 1}(z) \right] + \frac{\ell - L}{2} t_{p + 1}^2 \norm*{x^{p + 1} - y^{p + 1}}^2 + \left( a t_p - b + \frac{1}{4} \right) \sigma_p(z)
    .\end{multline}
    Now, let~$k \ge 0$.
    \Cref{thm:sigma k1 k2} with~$(k_1, k_2) = (p, k + 1)$ implies
    \begin{multline}
        t_p^2 \sigma_p(z) - t_{p + 1}^2 \sigma_{p + 1}(z) \ge \frac{\ell}{2} \left[ \rho_p(z) - \rho_{p - 1}(z) \right] \\
        + \frac{\ell - L}{2} t_{p + 1}^2 \norm*{x^{p + 1} - y^{p + 1}}^2 + \left( a t_p - b + \frac{1}{4} \right) \Biggl[ \sigma_{k + 1}(z) \\
        + \frac{\ell}{2} \left( \norm*{x^{k + 1} - x^k}^2 - \norm*{x^p - x^{p - 1}}^2 + \sum_{r = p}^k \frac{1}{t_r} \norm*{x^r - x^{r - 1}}^2 \right) \Biggr]
    .\end{multline}
    Adding up the above inequality from~$p = 1$ to~$p = k$, the fact that~$t_1 = 1$ and~$\rho_0(z) = \norm*{x^1 - z}^2$ leads to
    \begin{multline} \label{eq:key relation tmp}
        \sigma_1(z) - t_{k + 1}^2 \sigma_{k + 1}(z) \\
        \ge \frac{\ell}{2} \left[ \rho_{k}(z) - \norm*{x^1 - z}^2 \right] + \frac{\ell - L}{2} \sum_{p = 1}^k t_{k + 1}^2 \norm*{x^{k + 1} - y^{k + 1}}^2 \\
        + \left( a \sum_{p = 1}^{k} t_p + \left( \frac{1}{4} - b \right) k \right) \left[ \sigma_{k + 1}(z) + \frac{\ell}{2} \norm*{x^{k + 1} - x^k}^2 \right] \\
        - \frac{\ell}{2} \sum_{p = 1}^{k} \left( a t_p - b + \frac{1}{4} \right) \norm*{x^p - x^{p + 1}}^2 \\
        + \frac{\ell}{2} \sum_{p = 1}^{k} \left( a t_p - b + \frac{1}{4} \right) \sum_{r = p}^{k} \frac{1}{t_r} \norm*{x^r - x^{r - 1}}^2 
    .\end{multline}
    Let us write the last two terms of the right-hand side for~\cref{eq:key relation tmp} as~$S_1$ and~$S_2$, respectively.
    \Cref{eq:change sum} yields
    \[
        \begin{split}
            S_2 &= \frac{\ell}{2} \sum_{r = 1}^{k} \sum_{p = 1}^{r} \left( a t_p - b + \frac{1}{4} \right) \frac{1}{t_r} \norm*{x^r - x^{r - 1}}^2 \\
                &= \frac{\ell}{2} \sum_{p = 1}^{k} \sum_{r = 1}^{p} \left( a t_r - b + \frac{1}{4} \right) \frac{1}{t_p} \norm*{x^p - x^{p - 1}}^2
        .\end{split}
    \]
    Hence, it follows that
    \begin{multline} \label{eq:s_1 + s_2}
         S_1 + S_2 = \frac{\ell}{2} \sum_{p = 1}^{k} \left[ \frac{1}{t_p} \sum_{r = 1}^{p} \left( a t_r - b + \frac{1}{4} \right) - \left( a t_p - b + \frac{1}{4} \right) \right] \norm*{x^p - x^{p - 1}}^2 \\
         = \frac{\ell}{2} \sum_{p = 1}^{k} \frac{1}{t_p} \left[ a \left( \sum_{r = 1}^{p - 1} t_r - t_p^2 + t_p \right) + \left( \frac{1}{4} - b \right) (p - t_p) \right] \norm*{x^p - x^{p - 1}}^2
    .\end{multline}
    Again~$t_1 = 1$ gives
    \[
        \begin{split}
            - t_p^2 + t_p &= \sum_{r = 1}^{p - 1} ( - t_{r + 1}^2 + t_{r + 1} + t_r^2 - t_r ) = \sum_{r = 1}^{p - 1} \left(- (1 - a) t_r - b + \frac{1}{4} \right) \\
                          &= - (1 - a) \sum_{r = 1}^{p - 1} t_r + \left( \frac{1}{4} - b \right) (p - 1)
        ,\end{split}
    \]
    where the second equality comes from \cref{thm:t:t over-relax geq}.
    Thus, we get
    \[ \label{eq:sum t}
        \sum_{r = 1}^{p - 1} t_r = \frac{t_p^2 - t_p}{1 - a} + \left( \frac{1}{4} - b \right) \frac{p - 1}{1 - a}
    .\] 
    Substituting this into~\cref{eq:s_1 + s_2}, it follows that
    \begin{multline}
        S_1 + S_2 \\
        = \frac{\ell}{2 (1 - a)} \sum_{p = 1}^{k} \left[ a^2 (t_p - 1) + \left( \frac{1}{4} - b \right) \frac{p - t_p + a (t_p - 1)}{t_p} \right] \norm*{x^p - x^{p - 1}}^2
    .\end{multline}
    Combined with~\cref{eq:key relation tmp,eq:sum t}, we have
    \begin{align}
        \MoveEqLeft \sigma_1(z) - t_{k + 1}^2 \sigma_{k + 1}(z) \\
        \ge{}& \frac{\ell}{2} \left[ \rho_k(z) - \norm*{x^1 - z}^2 \right] + \frac{\ell - L}{2} \sum_{p = 1}^{k} t_{p + 1}^2 \norm*{x^{k + 1} - y^{k + 1}}^2 \\
             &+ \frac{1}{1 - a} \left[ a (t_{k + 1}^2 - t_{k + 1}) + \left( \frac{1}{4} - b \right) k \right] \left[ \sigma_{k + 1}(z) + \frac{\ell}{2} \norm*{x^{k + 1} - x^k}^2 \right] \\
             &+ \frac{\ell}{2 (1 - a)} \sum_{p = 1}^{k} \left[ a^2 (t_p - 1) + \left( \frac{1}{4} - b \right) \frac{p - t_p + a (t_p - 1)}{t_p} \right] \norm*{x^p - x^{p - 1}}^2
    .\end{align}
    Easy calculations give
    \begin{align}
        \MoveEqLeft \sigma_1(z) + \frac{\ell}{2} \norm*{x^1 - z}^2 \\
        \ge{}& \frac{1}{1 - a} \left[ t_{k + 1}^2 - a t_{k + 1} + \left( \frac{1}{4} - b \right) k \right] \sigma_{k + 1}(z) \\
        &+ \frac{\ell}{2 (1 - a)} \left[ a (t_{k + 1}^2 - t_{k + 1}) + \left( \frac{1}{4} - b \right) k \right] \norm*{x^{k + 1} - x^k}^2 \\
        &+ \frac{\ell}{2 (1 - a)} \sum_{p = 1}^{k} \left[ a^2 (t_p - 1) + \left( \frac{1}{4} - b \right) \frac{p - t_p + a(t_p - 1)}{t_p} \right] \norm*{x^p - x^{p - 1}}^2 \\
        &+ \frac{\ell}{2} \rho_k(z) + \frac{\ell - L}{2} \sum_{p = 1}^{k} t_{p + 1}^2 \norm*{x^{k + 1} - y^{k + 1}}^2
    .\end{align}
    \Cref{thm:sigma:1} with~$k = 0$ and~$y^1 = x^0$ and~\cref{eq:Pythagoras} with~$(v^1, v^2, v^3) = (x^0, x^1, z)$ lead to
    \[
        \sigma_1(z) \le - \frac{\ell}{2} \left[ \norm*{x^1 - z}^2 - \norm*{x^0 - z}^2 \right] - \frac{\ell - L}{2} \norm*{x^1 - y^1}^2
    .\]
    From the above two inequalities and the fact that~$\ell \ge L$, we can derive the desired inequality.
\end{proof}

Let us define the linear function~$P \colon \setR \to \setR$ and quadratic ones~$Q_1 \colon \setR \to \setR$,~$Q_2 \colon \setR \to \setR$, and~$Q_3 \colon \setR \to \setR$ by
\[ \label{eq:P Q}
    \begin{aligned}
    &P(\alpha) \coloneqq \frac{a^2 (\alpha - 1)}{2},\\
    &Q_1(\alpha) \coloneqq \frac{1 - a}{4} \alpha^2 + \left[ 1 - \frac{a}{2} + \frac{1 - 4 b}{4 (1 - a)} \right] \alpha + 1,\\
    &Q_2(\alpha) \coloneqq \frac{a (1 - a)}{4} \alpha^2 + \left[ \frac{a}{2} + \frac{1 - 4 b}{4 (1 - a)} \right] \alpha, \\
    &Q_3(\alpha) \coloneqq \left( \frac{1 - a}{2} \alpha + 1 \right)^2
    .\end{aligned}
\]
The following lemma provides the key relation to evaluate the convergence rate of \cref{alg:acc-pgm general}.
\begin{lemma} \label{thm:acc conv rate}
    Under \cref{asm:bound}, \cref{alg:acc-pgm general} generates a sequence~$\set*{x^k}$ such that
    \begin{multline}
        \frac{\ell R}{2} \ge Q_1(k) u_0(x^{k + 1}) + \frac{\ell}{2} Q_2(k) \norm*{x^{k + 1} - x^k}^2 + \frac{\ell}{2} \sum_{p = 1}^{k} P(p) \norm*{x^p - x^{p - 1}}^2 \\
         +  \frac{\ell - L}{2} \sum_{p = 1}^{k} Q_3(p) \norm*{x^{p + 1} - y^{p + 1}}^2
    \end{multline}
    for all~$k \ge 0$, where~$R \ge 0$ and~$P, Q_1, Q_2, Q_3 \colon \setR \to \setR$ are given in~\cref{eq:R,eq:P Q}, respectively, and~$u_0$ is a merit function defined by~\cref{eq:u_0}.
\end{lemma}
\begin{proof}    
    Let~$k \ge 0$.
    With similar arguments used in the proof of \cref{thm:conv rate} (see~\cite[Theorem 5.2]{Tanabe2022a}), we get
    \[
        \sup_{F^\ast \in F(X^\ast \cap \level_F(F(x^0)))} \inf_{z \in F^{-1}(\set*{F^\ast})} \sigma_{k + 1}(z) = u_0(x^{k + 1})
    .\] 
    Since~$\rho_k(z) \ge 0$, \cref{thm:key relation} and the above equality lead to
    \begin{align}
        \frac{\ell R}{2} \ge{}& \frac{1}{1 - a} \left[ t_{k + 1}^2 - a t_{k + 1} + \left( \frac{1}{4} - b \right) k \right] u_0(x^{k + 1}) \\
        &+ \frac{\ell}{2 (1 - a)} \left[ a (t_{k + 1}^2 - t_{k + 1}) + \left( \frac{1}{4} - b \right) k \right] \norm*{x^{k + 1} - x^k}^2 \\
        &+ \frac{\ell}{2 (1 - a)} \sum_{p = 1}^{k} \left[ a^2 (t_p - 1) + \left( \frac{1}{4} - b \right) \frac{p - t_p + a(t_p - 1)}{t_p} \right] \norm*{x^p - x^{p - 1}}^2 \\
        &+ \frac{\ell - L}{2} \sum_{p = 1}^{k} t_{p + 1}^2 \norm*{x^{p + 1} - y^{p + 1}}^2
    .\end{align}
    We now show that the coefficients of the four terms on the right-hand side can be bounded from below by the polynomials given in~\cref{eq:P Q}.
    First, by using the relation
    \[ \label{eq:t k+1 geq}
        t_{k + 1} \ge \frac{1 - a}{2} k + 1
    \] 
    obtained from \cref{thm:t:t geq} and~$a \in [0, 1)$, we have
    \begin{align}
        \MoveEqLeft \frac{1}{1 - a} \left[ t_{k + 1}^2 - a t_{k + 1} + \left( \frac{1}{4} - b \right) k \right] = \frac{1}{1 - a} \left[ t_{k + 1} (t_{k + 1} - a) + \left( \frac{1}{4} - b \right) k \right] \\
        &\ge \frac{1}{1 - a} \left[ \left( \frac{1 - a}{2} k + 1 \right) \left( \frac{1 - a}{2} k + 1 - a \right) + \left( \frac{1}{4} - b \right) k \right] = Q_1(k)
    .\end{align}
    Again,~\cref{eq:t k+1 geq} gives
    \begin{align}
        \MoveEqLeft \frac{1}{1 - a} \left[ a (t_{k + 1}^2 - t_{k + 1}) + \left( \frac{1}{4} - b \right) k \right] 
        = \frac{a}{1 - a} t_{k + 1} (t_{k + 1} - 1) + \frac{1 - 4 b}{4 (1 - a)} k \\
        &\ge \frac{a}{1 - a} \left( \frac{1 - a}{2} k + 1 \right) \left( \frac{1 - a}{2} k \right) + \frac{1 - 4 b}{4 (1 - a)} k = Q_2(k)
    .\end{align}
    Moreover, since~$t_p \le p$ (cf.~\cref{thm:t:t leq}),~$t_k \ge 1$ (cf.~\cref{thm:t:t geq}), and~$b \in (a^2 / 4, 1 / 4]$, we obtain
    \[
        \frac{1}{1 - a} \left[ a^2 (t_p - 1) + \left( \frac{1}{4} - b \right) \frac{p - t_p + a (t_p - 1)}{t_p} \right] \ge \frac{a^2}{1 - a} (t_p - 1) \ge P(p)
    .\] 
    It is also clear from~\cref{eq:t k+1 geq} that
    \[
        t_{p + 1}^2 \ge Q_3(p)
    .\] 
    Thus, combining the above five inequalities, we get the desired inequality.
\end{proof}

Then, we can finally prove the main theorem.
\begin{proof}[\cref{thm:main theomrem of stepsize section:O(1/k2)}]
    It is clear from \cref{thm:acc conv rate} and~$Q_1(k) = O(k^2)$ as~$k \to \infty$.
\end{proof}
\begin{remark}
    \Cref{thm:acc conv rate} also implies the following other claims than \cref{thm:main theomrem of stepsize section:O(1/k2)}:
    \begin{itemize}
        \item $O(1 / k^2)$ convergence rate of~$\set*{\norm*{x^{k + 1} - x^k}^2}$ when~$a > 0$;
        \item the absolute convergence of~$\set*{k \norm*{x^{k + 1} - x^k}^2}$ when~$a > 0$;
        \item the absolute convergence of~$\set*{k^2 \norm*{x^k - y^k}^2}$ when~$\ell > L$.
    \end{itemize}
    Note that the second one generalize~\cite[Corollary~3.2]{Chambolle2015} for single-objective problems. 
\end{remark}

\ifthenelse{\boolean{isMain}}{ 
}{
   {\bibliographystyle{jorsj} 
   \bibliography{library}} 
}

\section{Convergence of the iterates} \label{sec: convergence sequence}
While the last section shows that \cref{alg:acc-pgm general} has an~$O(1 / k^2)$ convergence rate like \cref{alg:acc-pgm}, this section proves the following theorem, which is more strict than \cref{thm:accumulation point} related to \cref{alg:acc-pgm}:
\begin{theorem} \label{thm:main convergence}
    Let~$\set*{x^k}$ be generated by \cref{alg:acc-pgm general} with~$a > 0$.
    Then, under \cref{asm:bound}, the following two properties hold:
    \begin{enumerate}
        \item \label{thm:main convergence:bound} $\set*{x^k}$ is bounded, and it has an accumulation point;
        \item \label{thm:main convergence:Pareto} $\set*{x^k}$ converges to a weak Pareto optimum for~\cref{eq:MOP}.
    \end{enumerate}
\end{theorem}
The latter claim is also significant in application.
For example, finite-time manifold (active set) identification, which detects the low-dimensional manifold where the optimal solution belongs, essentially requires only the convergence of the generated sequence to a unique point rather than the strong convexity of the objective functions~\citep{Sun2019}.

Again, we will prove \cref{thm:main convergence} after showing some lemmas.
First, we mention the following result, obvious from \cref{asm:bound,thm:main theorem of stepsize section:less than x0}.
\begin{lemma} \label{thm:z exist}
    Let~$\set*{x^k}$ be generated by \cref{alg:acc-pgm general}.
    Then, for any~$k \ge 0$, there exists~$z \in X^\ast \cap \level_F(F(x^0))$ (see \cref{eq:weak Pareto,eq:level set} for the definitions of~$X^\ast$ and~$\level_F$) such that
    \[
        \sigma_k(z) \ge 0 \eqand \norm*{z - x^0}^2 \le R
    ,\] 
    where~$R \ge 0$ is given by~\cref{eq:R}.
\end{lemma}

The following lemma also contributes strongly to the proof of the main theorem.
\begin{lemma} \label{thm:prod}
    Let~$\set*{\gamma_q}$ be defined by \cref{line:gamma} in \cref{alg:acc-pgm general}.
    Then, we have
    \[
        \sum_{p = s}^r \prod_{q = s}^p \gamma_q \le 2 ( s - 1) \forallcondition{s, r \ge 1}
    .\] 
\end{lemma}
\begin{proof}
    By using \Cref{thm:t:gamma}, we see that
    \[
        \prod_{q = s}^{p} \gamma_q \le \prod_{q = s}^{p} \frac{q - 1}{q + 1 / 2} 
    .\]
    Let~$\Gamma$ and~$\Beta$ denote the gamma and beta functions defined by
    \[ \label{eq:gamma and beta}
        \Gamma(\alpha) \coloneqq \int_{0}^{\infty} \tau^{\alpha - 1} \exp(-\tau) \odif{\tau} \eqand
        \Beta(\alpha, \beta) \coloneqq \int_{0}^{1} \tau^{\alpha - 1} (1 - \tau)^{\beta - 1} \odif{\tau}
    ,\]
    respectively.
    Applying the well-known properties:
    \[ \label{eq:gamma and beta properties}
        \Gamma(\alpha) = (\alpha - 1)!, \quad \Gamma(\alpha + 1) = \alpha \Gamma(\alpha), \eqand B(\alpha, \beta) = \frac{\Gamma(\alpha) \Gamma(\beta)}{\Gamma(\alpha + \beta)}
    .\]
    we get
    \[
        \prod_{q = s}^{p} \gamma_q \le \frac{\Gamma(p) / \Gamma(s - 1)}{\Gamma(p + 3 / 2) / \Gamma(s + 1 / 2)}
        = \frac{B(p, 3 / 2)}{B(s - 1, 3 / 2)}
    .\] 
    This implies
    \[
        \sum_{p = s}^{r} \prod_{q = s}^{p} \gamma_q \le \sum_{p = 1}^{r} B(p, 3 / 2) / B(s - 1, 3 / 2)
    .\] 
    Then, it follows from the definition~\cref{eq:gamma and beta} of~$\Beta$ that
    \[
        \begin{split}
            \sum_{p = s}^{r} \prod_{q = s}^{p} \gamma_q 
        &\le \sum_{p = s}^{r} \int_{0}^{1} \tau^{p - 1} (1 - \tau)^{1 / 2} \odif{\tau} / B(s - 1, 3 / 2) \\
        &= \int_{0}^{1} \sum_{p = s}^{r} \tau^{p - 1} (1 - \tau)^{1 / 2} \odif{\tau} / B(s - 1, 3 / 2) \\
        &= \int_{0}^{1} \frac{\tau^{s - 1} - \tau^r}{1 - \tau} (1 - \tau)^{1 / 2} \odif{\tau} / B(s - 1, 3 / 2) \\
        &= \frac{B(s, 1 / 2) - B(r + 1, 1 / 2)}{B(s - 1, 3 / 2)} 
        \le \frac{B(s, 1 / 2)}{B(s - 1, 3 / 2)}
        .\end{split}
    \] 
    Using again~\cref{eq:gamma and beta properties}, we conclude that
    \[
        \sum_{p = s}^{r} \prod_{q = s}^{p} \gamma_q
        \le \frac{\Gamma(s) \Gamma(1 / 2) / \Gamma(s + 1 / 2)}{\Gamma(s - 1) \Gamma(3 / 2) / \Gamma(s + 1 / 2)}
        = 2 (s - 1)
    .\] 
\end{proof}

Now, we introduce two functions~$\omega_k \colon \setR^n \to \setR$ and~$\nu_k \colon \setR^n \to \setR$ for any~$k \ge 1$, which will help our analysis, by
\begin{align}
    \omega_k(z) &\coloneqq \max \left( 0, \norm*{x^k - z}^2 - \norm*{x^{k - 1} - z}^2 \right) \label{eq:omega}, \\
    \nu_k(z) &\coloneqq \norm*{x^k - z}^2 - \sum_{s = 1}^{k} \omega_s(z) \label{eq:nu}
.\end{align}
The lemma below describes the properties of~$\omega_k$ and~$\nu_k$.
\begin{lemma} \label{thm:omega nu}
    Let~$\set*{x^k}$ be generated by \cref{alg:acc-pgm general} and recall that~$X^\ast, \level_F, \omega_k$, and~$\nu_k$ are defined by~\cref{eq:weak Pareto,eq:level set,eq:omega,eq:nu}, respectively.
    Moreover, suppose that \cref{asm:bound} holds and that~$z \in X^\ast \cap \level_F(F(x^0))$ satisfies the statement of \cref{thm:z exist} for some~$k \ge 1$.
    Then, it follows for all~$r = 1, \dots, k$ that
    \begin{enumerate}
        \item $\displaystyle \sum_{s = 1}^{r} \omega_s(z) \le \sum_{s = 1}^{r}(6s - 5) \norm*{x^s - x^{s - 1}}^2;$ \label{thm:omega nu:omega}
        \item $\displaystyle \nu_{r + 1}(z) \le \nu_r(z).$ \label{thm:omega nu:nu}
    \end{enumerate}
\end{lemma}
\begin{proof}
    \sublabelcref{thm:omega nu:omega}:
    Let~$k \ge p \ge 1$.
    From the definition of~$y^{p + 1}$ given in \cref{line:y} of \cref{alg:acc-pgm}, we have
    \[
        \begin{split}
            &\norm*{x^{p + 1} - z}^2 - \norm*{x^p - z}^2 \\
        &= - \norm*{x^{p + 1} - x^p}^2 + 2 \innerp*{x^{p + 1} - y^{p + 1}}{x^{p + 1} - z} + 2 \gamma_p \innerp*{x^p - x^{p - 1}}{x^{p + 1} - z} \\
        &= - \norm*{x^{p + 1} - x^p}^2 + 2 \innerp*{x^{p + 1} - y^{p + 1}}{y^{p + 1} - z} + 2 \norm*{x^{p + 1} - y^{p + 1}}^2 \\
        &\quad + 2 \gamma_p \innerp*{x^p - x^{p - 1}}{x^{p + 1} - z}
     .   \end{split}
    \]
    On the other hand, \cref{thm:sigma:1} gives
    \[
        2 \innerp*{x^{p + 1} - y^{p + 1}}{y^{p + 1} - z} \le - \frac{2}{\ell} \sigma_{p + 1}(z) - \frac{2 \ell - L}{\ell} \norm*{x^{p + 1} - y^{p + 1}}^2
    .\]
    Moreover, \cref{thm:sigma k1 k2} with~$(k_1, k_2) = (p + 1, k + 1)$ implies
    \[
        \begin{split}
            \MoveEqLeft
        - \frac{2}{\ell} \sigma_{p + 1}(z) \\
        &\le - \frac{2}{\ell} \sigma_{k + 1}(z) - \norm*{x^{k + 1} - x^k}^2 + \norm*{x^{p + 1} - x^p}^2 - \sum_{r = p + 1}^{k} \frac{1}{t_r}\norm*{x^r - x^{r - 1}}^2 \\
        &\le \norm*{x^{p + 1} - x^p}^2
        ,\end{split}
    \]
    where the second inequality comes from the assumption on~$z$.
    Combining the above three inequalities, we get
    \begin{multline}
        \norm*{x^{p + 1} - z}^2 - \norm*{x^p - z}^2 \le \frac{L}{\ell} \norm*{x^{p + 1} - y^{p + 1}}^2 + 2 \gamma_p \innerp*{x^p - x^{p - 1}}{x^{p + 1} - z} \\
        \shoveleft = \frac{L}{\ell} \norm*{x^{p + 1} - y^{p + 1}}^2
        + \gamma_p \Bigl( \norm*{x^p - z}^2 - \norm*{x^{p - 1} - z}^2 + \norm*{x^p - x^{p - 1}}^2 \\
    + 2 \innerp*{x^p - x^{p - 1}}{x^{p + 1} - x^p} \Bigr)
    .\end{multline}
    Using the relation~$\norm*{x^{p + 1} - y^{p + 1}}^2 + 2 \gamma_p \innerp*{x^p - x^{p - 1}}{x^{p + 1} - x^p} = \norm*{x^{p + 1} - x^p}^2 + \gamma_p^2 \norm*{x^p - x^{p - 1}}^2$, which holds from the definition of~$y^k$, we have
    \begin{multline}
        \norm*{x^{p + 1} - z}^2 - \norm*{x^p - z}^2
        \le - \frac{\ell - L}{\ell} \norm*{x^{p + 1} - y^{p + 1}}^2 + \norm*{x^{p + 1} - x^p}^2 \\
        + \gamma_p \left( \norm*{x^p - z}^2 - \norm*{x^{p - 1} - z}^2 \right) + ( \gamma_p + \gamma_p^2 ) \norm*{x^p - x^{p - 1}}^2
    .\end{multline}
    Since~$0 \le \gamma_p \le 1$ from \cref{thm:t:gamma} and~$\ell \ge L$, we obtain
    \begin{multline}
        \norm*{x^{p + 1} - z}^2 - \norm*{x^p - z}^2 \\
        \le \gamma_p \left( \norm*{x^p - z}^2 - \norm*{x^{p - 1} - z}^2 + 2 \norm*{x^p - x^{p - 1}}^2 \right) + \norm*{x^{p + 1} - x^p}^2 \\
        \le \gamma_p \left( \omega_p(z) + 2 \norm*{x^p - x^{p - 1}}^2 \right) + \norm*{x^{p + 1} - x^p}^2
    ,\end{multline}
    where the second inequality follows from the definition~\cref{eq:omega} of~$\omega_p$.
    Since the right-hand side is nonnegative,~\cref{eq:omega} again gives
    \[
        \omega_{p + 1}(z) \le \gamma_p \left( \omega_p(z) + 2 \norm*{x^p - x^{p - 1}}^2 \right) + \norm*{x^{p + 1} - x^p}^2 
    .\]
    Let~$s \le k$.
    Applying the above inequality recursively and using~$\gamma_1 = 0$, we get
    \[
        \begin{split}
            \omega_s(z) &\le 3 \sum_{p = 2}^{s} \prod_{q = p}^{s} \gamma_q \norm*{x^p - x^{p - 1}}^2 + 2 \prod_{q = 1}^{s} \gamma_q \norm*{x^1 - x^0}^2 + \norm*{x^s - x^{s - 1}}^2 \\
                        &\le 3 \sum_{p = 2}^{s} \prod_{q = p}^{s} \gamma_q \norm*{x^p - x^{p - 1}}^2 + \norm*{x^s - x^{s - 1}}^2
        .\end{split}
    \]
    Adding up the above inequality from~$s = 1$ to~$s = r \le k$, we have
    \[
        \begin{split}
        \sum_{s = 1}^r \omega_s(z)
        &\le 3 \sum_{s = 1}^r \sum_{p = 1}^s \prod_{q = p}^s \gamma_q \norm*{x^p - x^{p - 1}}^2 + \sum_{s = 1}^r \norm*{x^s - x^{s - 1}}^2 \\
        &= 3 \sum_{p = 1}^r \sum_{s = p}^r \prod_{q = p}^s \gamma_q \norm*{x^p - x^{p - 1}}^2 + \sum_{s = 1}^r \norm*{x^s - x^{s - 1}}^2 \\
        &= \sum_{s = 1}^r \left( 3 \sum_{p = s}^r \prod_{q = s}^p \gamma_q + 1 \right) \norm*{x^s - x^{s - 1}}^2
        ,\end{split}
    \]
    where the first equality follows from~\cref{eq:change sum}.
    Thus, \cref{thm:prod} implies
    \[
        \sum_{s = 1}^{r} \omega_s(z) \le \sum_{s = 1}^{r} (6 s - 5) \norm*{x^s - x^{s - 1}}^2 
    .\]

    \sublabelcref{thm:omega nu:nu}:
    \Cref{eq:nu} yields
    \[
        \begin{split}
            \nu_{r + 1}(z) &= \norm*{x^{r + 1} - z}^2 - \omega_{r + 1}(z) - \sum_{s = 1}^r \omega_s(z) \\
                           &= \norm*{x^{r + 1} - z}^2 - \max \left( 0, \norm*{x^{r + 1} - z}^2 - \norm*{x^r - z}^2 \right) - \sum_{s = 1}^{r} \omega_s(z) \\
                           &\le \norm*{x^{r + 1} - z}^2 - \left( \norm*{x^{r + 1} - z}^2 - \norm*{x^r - z}^2 \right) - \sum_{s = 1}^{r} \omega_s(z) \\
                           &= \norm*{x^r - z}^2 - \sum_{s = 1}^{r} \omega_s(z) = \nu_r(z) 
        ,\end{split}
    \] 
    where the second and third equalities come from the definitions~\cref{eq:omega,eq:nu} of~$\omega_{r + 1}$ and~$\nu_r$, respectively.
\end{proof}

Let us now prove the first part of the main theorem.
\begin{proof}[\cref{thm:main convergence:bound}]
    Let~$k \ge 1$ and suppose that~$z \in X^\ast \cap \level_F(F(x^0))$ satisfies the statement of \cref{thm:z exist}, where~$X^\ast$ and~$\level_F$ are given by~\cref{eq:weak Pareto,eq:level set}, respectively.
    Then, \cref{thm:omega nu:nu} gives
    \[
        \begin{split}
            \nu_k(z) &\le \nu_1(z) = \norm*{x^1 - z}^2 - \omega_1(z) \\
                           &= \norm*{x^1 - z}^2 - \max \left( 0, \norm*{x^1 - z}^2 - \norm*{x^0 - z}^2 \right) \\
                           &\le \norm*{x^1 - z}^2 - \left( \norm*{x^1 - z}^2 - \norm*{x^0 - z}^2 \right) = \norm*{x^0 - z}^2
        ,\end{split}
    \] 
    where the second equality follows from the definition~\cref{eq:omega} of~$\omega_1$.
    Considering the definition~\cref{eq:nu} of~$\nu_k$, we obtain
    \[
        \norm*{x^k - z}^2 \le \norm*{x^0 - z}^2 + \sum_{s = 1}^{k} \omega_s(z)
    .\] 
    Taking the square root of both sides and using~\cref{eq:omega}, we get
    \[
        \norm*{x^k - z} \le \sqrt{ \norm*{x^0 - z}^2 + \sum_{s = 1}^{k} (6s - 5) \norm*{x^s - x^{s - 1}}^2 }
    .\] 
    Applying the reverse triangle inequality~$\norm*{x^k - x^0} - \norm*{x^0 - z} \le \norm*{x^k - z}$ to the left-hand side leads to
    \begin{align}
        \norm*{x^k - x^0} &\le \norm*{x^0 - z} + \sqrt{ \norm*{x^0 - z}^2 + \sum_{s = 1}^{k} (6s - 5) \norm*{x^s - x^{s - 1}}^2 } \\
                          &\le \sqrt{R} + \sqrt{R + \sum_{s = 1}^{k} (6s - 5) \norm*{x^s - x^{s - 1}}^2}
    ,\end{align}
    where the second inequality comes from the assumption on~$z$.
    Moreover, since~$a > 0$, the right-hand side is bounded from above according to \cref{thm:acc conv rate}.
    This implies that~$\set*{x^k}$ is bounded, and so it has accumulation points.
\end{proof}

Before proving \cref{thm:main convergence:Pareto}, we show the following lemma.
\begin{lemma} \label{thm:convergence norm}
    Let~$\set*{x^k}$ be generated by \cref{alg:acc-pgm general} with~$a > 0$ and suppose that \cref{asm:bound} holds.
    Then, if~$\bar{z}$ is an accumulation point of~$\set*{x^k}$, then~$\set*{\norm*{x^k - \bar{z}}}$ is convergent.
\end{lemma}
\begin{proof}
    Assume that~$\set*{x^{k_j}} \subseteq \set*{x^k}$ converges to~$\bar{z}$.
    Then, we have~$\sigma_{k_j}(\bar{z}) \to 0$ by the definition~\cref{eq:sigma rho} of~$\sigma_{k_j}$.
    Therefore, we can regard~$\bar{z}$ to satisfy the statement of \cref{thm:z exist} at~$k = \infty$, and thus the inequalities of \cref{thm:omega nu} hold for any~$r \ge 1$ and~$\bar{z}$.
    This means~$\set*{\nu_k(\bar{z})}$ is non-increasing and bounded, i.e., convergent.
    Hence~$\set*{\norm*{x^k - \bar{z}}}$ is convergent.
\end{proof}

Finally, we finish the proof of the main theorem.
\begin{proof}[\cref{thm:main convergence:Pareto}]
    Suppose that~$\set*{x^{k^1_j}}$ and~$\set*{x^{k^2_j}}$ converges to~$\bar{z}^1$ and~$\bar{z}^2$, respectively.
    From \cref{thm:convergence norm}, we see that
    \[
        \lim_{j \to \infty} \left( \norm*{x^{k^2_j} - \bar{z}^1}^2 - \norm*{x^{k^2_j} - \bar{z}^2}^2 \right) = \lim_{j \to \infty} \left( \norm*{x^{k^1_j} - \bar{z}^1}^2 - \norm*{x^{k^1_j} - \bar{z}^2}^2 \right) 
    .\] 
    This yields that~$\norm*{\bar{z}^1 - \bar{z}^2}^2 = - \norm*{\bar{z}^1 - \bar{z}^2}^2$, and so~$\norm*{\bar{z}^1 - \bar{z}^2}^2 = 0$, i.e.,~$\set*{x^k}$ is convergent.
    Let~$x^k \to x^\ast$.
    Since~$\norm*{x^{k + 1} - x^k}^2 \to 0$, $\set*{y^k}$ is also convergent to~$x^\ast$.
    Therefore, \cref{thm:acc prox termination} shows that~$x^\ast$ is weakly Pareto optimal for~\cref{eq:MOP}.
\end{proof}

\ifthenelse{\boolean{isMain}}{ 
}{
   {\bibliographystyle{jorsj} 
   \bibliography{library}} 
}

\section{Numerical experiments} \label{sec: experiments}
This section compares the performance between \cref{alg:acc-pgm general} with various~$a$ and~$b$ and \cref{alg:acc-pgm} ($a = 0, b = 1 / 4$) through numerical experiments.
We run all experiments in Python 3.9.9 on a machine with 2.3 GHz Intel Core i7 CPU and 32 GB memory.
For each example, we test 15 different hyperparameters combining~$a = 0, 1 / 6, 1 / 4, 1 / 2, 3 / 4$ and~$b = a^2 / 4, (a^2 + 1) / 8, 1 / 4$, i.e.,
\[
    (a, b) = \left\{
        \begin{gathered}
            (0, 0), (0, 1 / 8), (0, 1 / 4),\\
            (1 / 6, 1 / 144), (1 / 6, 37 / 288), (1 / 6, 1 / 4),\\
            (1 / 4, 1 / 64), (1 / 4, 17 / 128), (1 / 4, 1 / 4), \\
            (1 / 2, 1 / 16), (1 / 2, 5 / 32), (1 / 2, 1 / 4), \\
            (3 / 4, 9 / 64), (3 / 4, 25 / 128), (3 / 4, 1 / 4)
        \end{gathered}
    \right\},
\] 
and we set~$\varepsilon = 10^{-5}$ for the stopping criteria.

\subsection{Artificial test problems (bi-objective and tri-objective)}
First, we solve the multi-objective test problems in the form~\cref{eq:MOP} used in~\cite{Tanabe2022a}, modifications from~\cite{Jin2001,Fliege2009}, whose objective functions are defined by
\begin{align}
    &f_1(x) = \frac{1}{n} \norm*{x}^2,
    f_2(x) = \frac{1}{n} \norm*{x - 2}^2,
    g_1(x) = g_2(x) = 0, \tag{JOS1}\label{eq:JOS1} \\
    &f_1(x) = \frac{1}{n} \norm*{x}^2,
    f_2(x) = \frac{1}{n} \norm*{x - 2}^2,
    g_1(x) = \frac{1}{n} \norm*{x}_1,
    g_2(x) = \frac{1}{2n} \norm*{x - 1}_1,
    \label{eq:JOS1_L1} \tag{JOS1-L1}\\
    &\left\{\begin{aligned} 
            f_1(x) &= \frac{1}{n^2} \sum_{i = 1}^{n} i (x_i - i)^4,
            f_2(x) = \exp \left( \sum_{i = 1}^{n} \frac{x_i}{n} \right) + \norm*{x}^2,\\
            f_3(x) &= \frac{1}{n (n + 1)} \sum_{i = 1}^{n} i (n - i + 1) \exp (- x_i),
            g_1(x) = g_2(x) = g_3(x) = 0,
    \end{aligned} \right. \label{eq:FDS} \tag{FDS}\\
    &\left\{\begin{aligned} 
            f_1(x) &= \frac{1}{n^2} \sum_{i = 1}^{n} i (x_i - i)^4,
            f_2(x) = \exp \left( \sum_{i = 1}^{n} \frac{x_i}{n} \right) + \norm*{x}^2,\\
            f_3(x) &= \frac{1}{n (n + 1)} \sum_{i = 1}^{n} i (n - i + 1) \exp (- x_i),
            g_1(x) = g_2(x) = g_3(x) = \indicator_{\setRpos^n}(x),
    \end{aligned} \right.\label{eq:FDS_CONSTRAINED} \tag{FDS-CON}
\end{align}
where~$x \in \setR^n, n = 50$ and~$\indicator_{\setRpos^n}$ is an indicator function~\cref{eq:indicator} of the nonnegative orthant.
We choose~$1000$ initial points, commonly for all pairs~$(a, b)$, and randomly with a uniform distribution between~$\underline{c}$ and~$\overline{c}$, where~$\underline{c} = (-2, \dots, -2)^\T$ and~$\overline{c} = (4, \dots, 4)^\T$ for~\cref{eq:JOS1,eq:JOS1_L1},~$\underline{c} = (-2, \dots, -2)^\T$ and~$\overline{c} = (2, \dots, 2)^\T$ for~\cref{eq:FDS}, and~$\underline{c} = (0, \dots, 0)^\T$ and~$\overline{c} = (2, \dots, 2)^\T$ for~\cref{eq:FDS_CONSTRAINED}.
Moreover, we use backtracking for updating~$\ell$, with~$1$ as the initial value of~$\ell$ and~$2$ as the constant multiplied into~$\ell$ at each iteration (cf.~\cite[Remark~4.1~(v)]{Tanabe2022a}).
Furthermore, at each iteration, we transform the subproblem~\cref{eq:acc prox subprob} into their dual as suggested in~\cite{Tanabe2022a} and solve them with the trust-region interior point method~\citep{Byrd1999} using the scientific library SciPy.

\begin{figure}[htbp]
    \centering
    \begin{minipage}[b]{.49\hsize}
        \centering
        \begin{minipage}[b]{.49\hsize}
            \centering
            \adjincludegraphics[trim={{.66\width} {.8\height} 0 0}, clip, width=\linewidth]{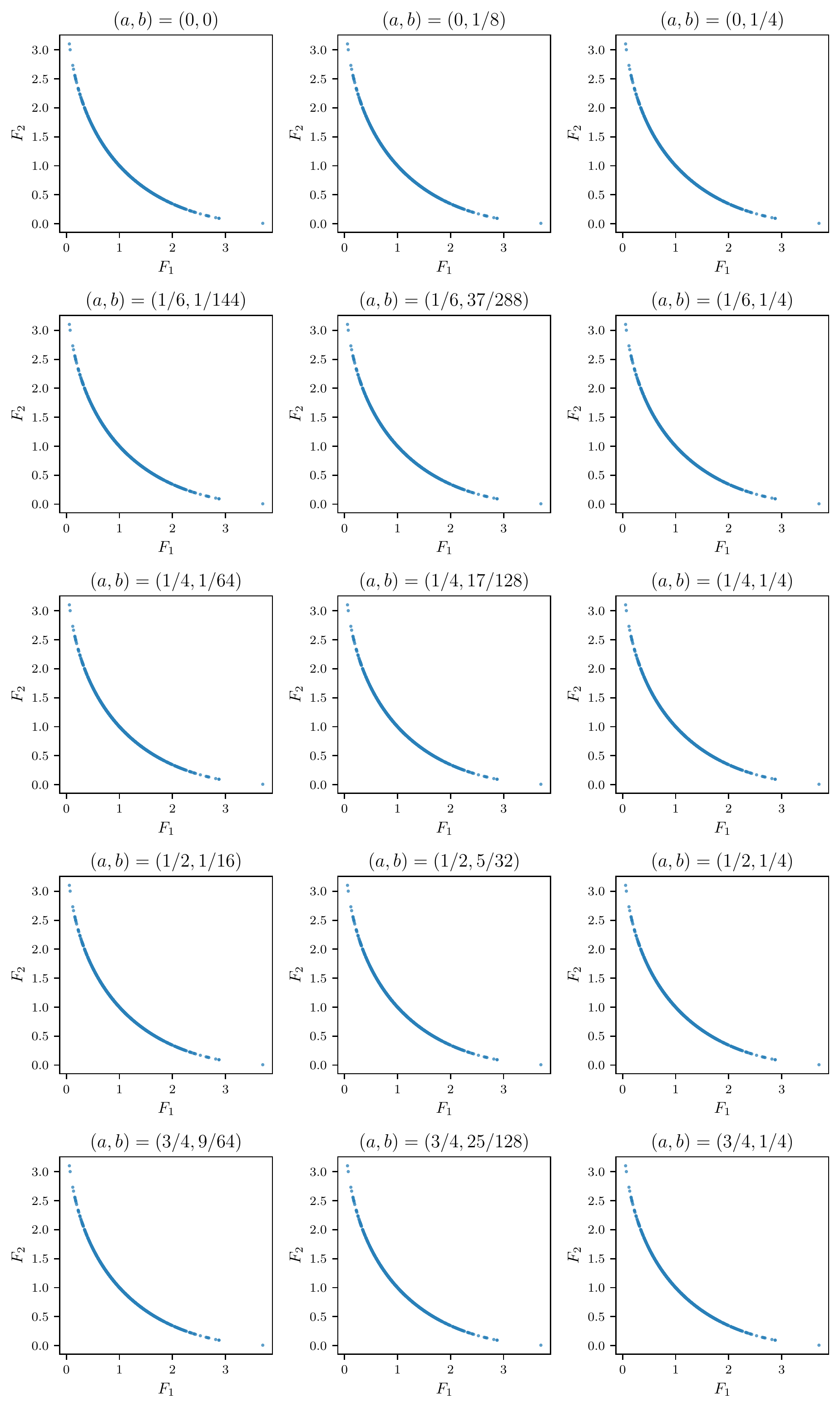}
        \end{minipage}
        \begin{minipage}[b]{.49\hsize}
            \centering
            \adjincludegraphics[trim={{.66\width} 0 0 {.8\height}}, clip, width=\linewidth]{figs/JOS1_ab.pdf}
        \end{minipage}
        \subcaption{\cref{eq:JOS1}}
        \label{fig:JOS1}
    \end{minipage}
    \begin{minipage}[b]{.49\hsize}
        \centering
        \begin{minipage}[b]{.49\hsize}
            \centering
            \adjincludegraphics[trim={{.66\width} {.8\height} 0 0}, clip, width=\linewidth]{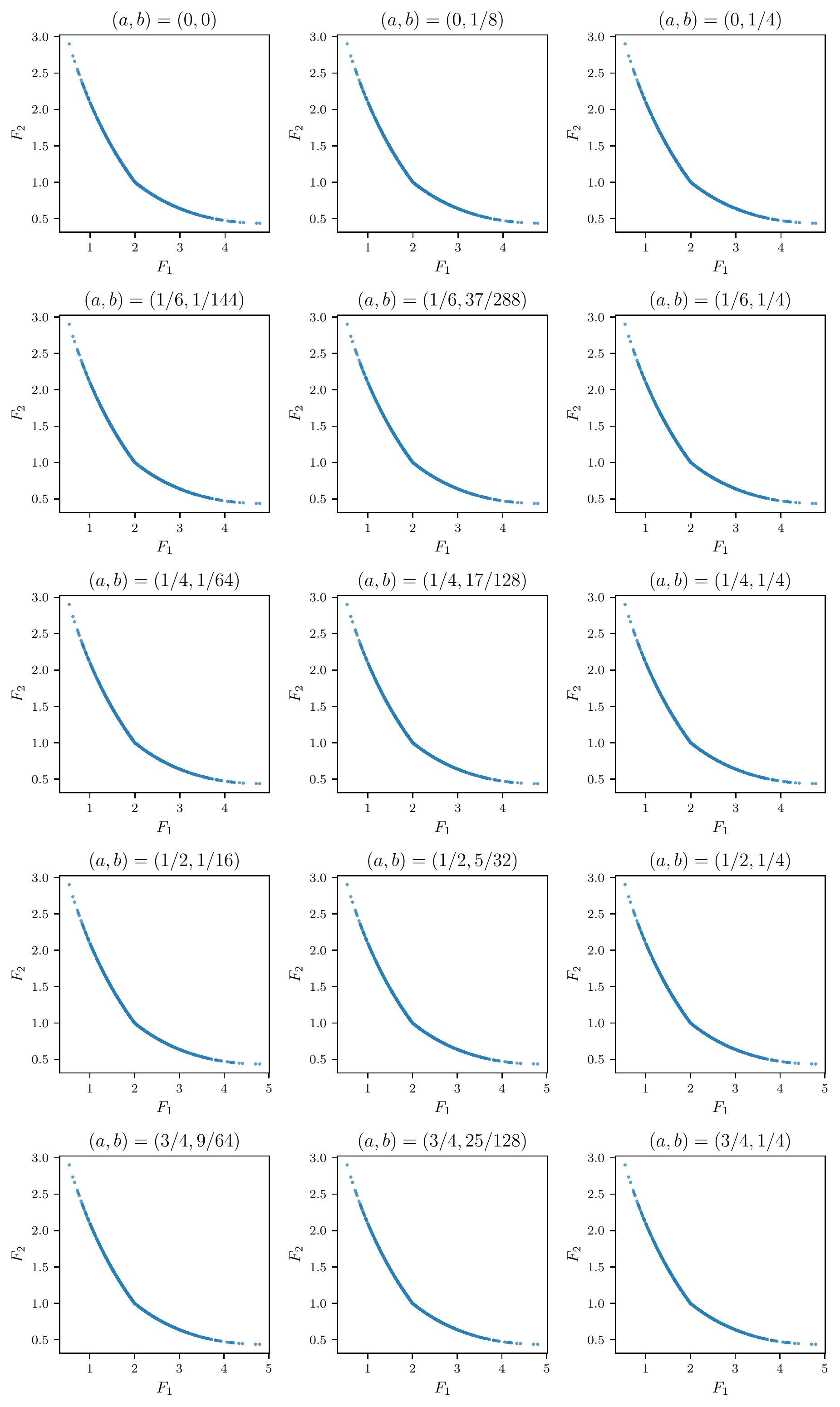}
        \end{minipage}
        \begin{minipage}[b]{.49\hsize}
            \centering
            \adjincludegraphics[trim={{.66\width} 0 0 {.8\height}}, clip, width=\linewidth]{figs/JOS1_L1_ab.pdf}
        \end{minipage}
        \subcaption{\cref{eq:JOS1_L1}}
        \label{fig:JOS1_L1}
    \end{minipage}
    \begin{minipage}[b]{.49\hsize}
        \centering
        \begin{minipage}[b]{.49\hsize}
            \centering
            \adjincludegraphics[trim={{.605\width} {.8\height} {.014\width} 0}, clip, width=\linewidth]{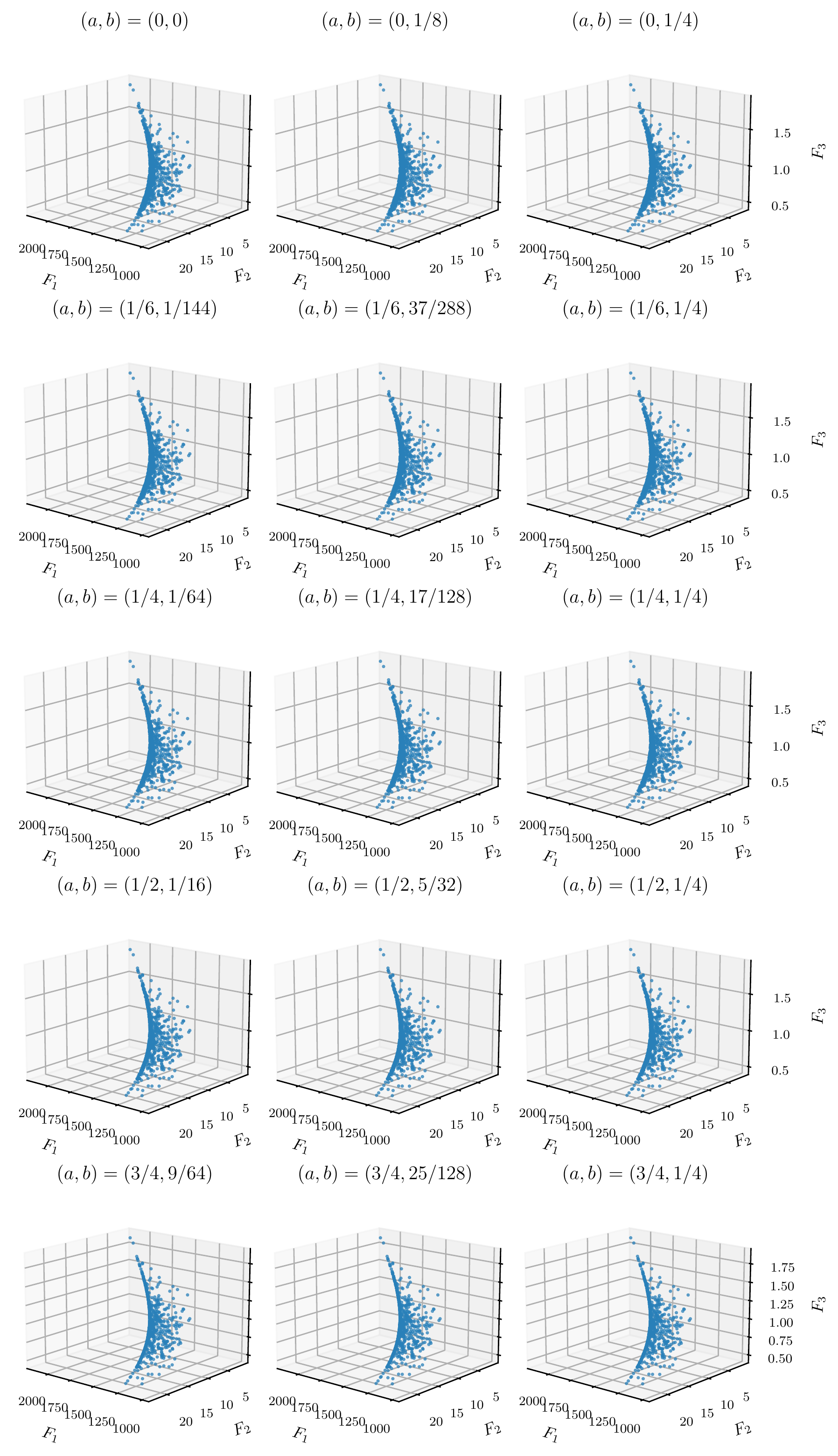}
        \end{minipage}
        \begin{minipage}[b]{.49\hsize}
            \centering
            \adjincludegraphics[trim={{.605\width} 0 {.014\width} {.8\height}}, clip, width=\linewidth]{figs/fDS_ab.pdf}
        \end{minipage}
        \subcaption{\cref{eq:FDS}}
        \label{fig:FDS}
    \end{minipage}
    \begin{minipage}[b]{.49\hsize}
        \centering
        \begin{minipage}[b]{.49\hsize}
            \centering
            \adjincludegraphics[trim={{.605\width} {.8\height} {.014\width} 0}, clip, width=\linewidth]{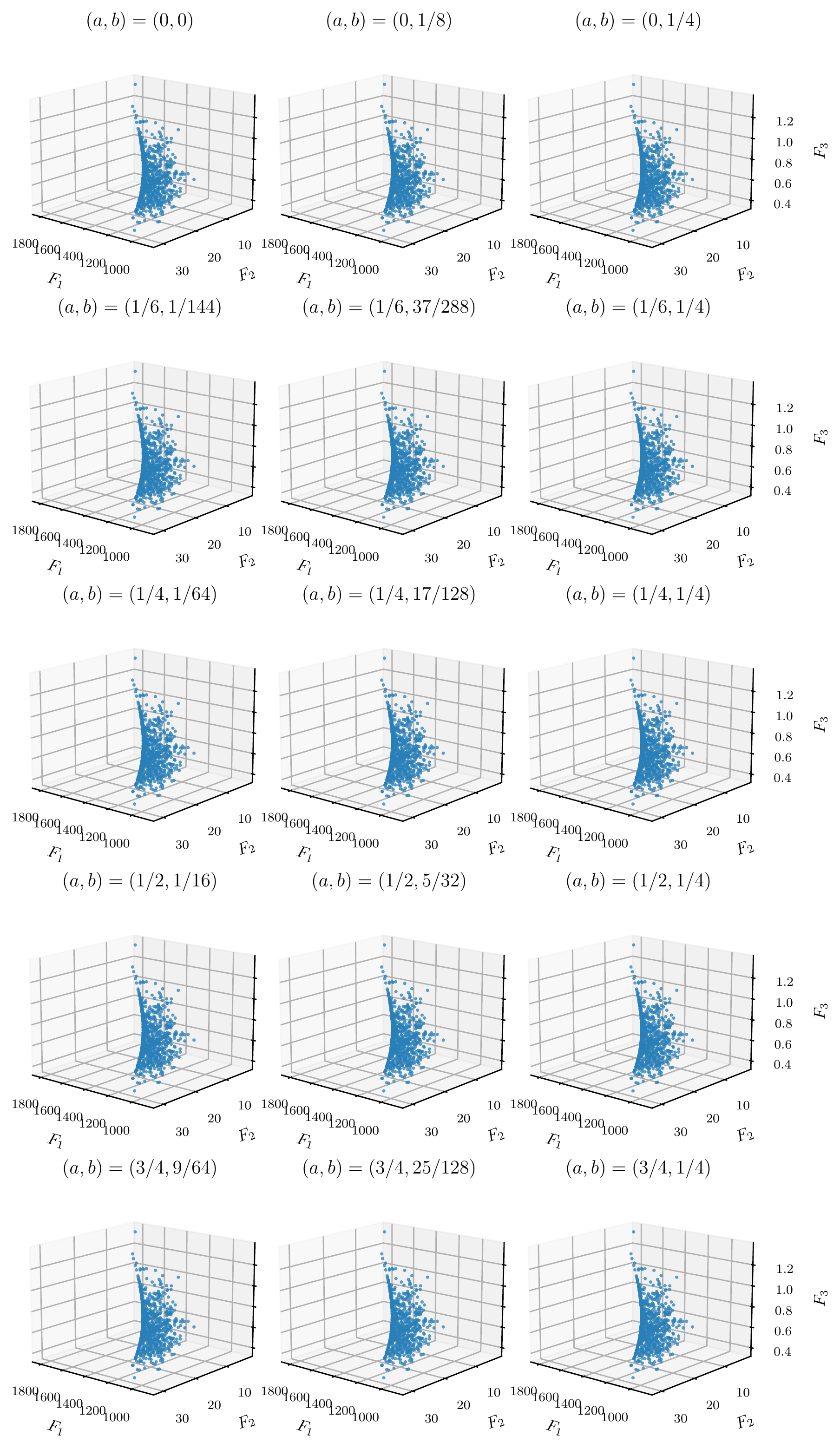}
        \end{minipage}
        \begin{minipage}[b]{.49\hsize}
            \centering
            \adjincludegraphics[trim={{.605\width} 0 {.014\width} {.8\height}}, clip, width=\linewidth]{figs/FDS_CONSTRAINED_ab.pdf}
        \end{minipage}
        \subcaption{\cref{eq:FDS_CONSTRAINED}}
        \label{fig:FDS_CON}
    \end{minipage}
    \caption{Pareto solutions obtained with some~$(a, b)$}
    \label{fig:Pareto}
\end{figure}

\Cref{fig:Pareto,tab:Average computational costs} present the experimental results.
\Cref{fig:Pareto} plots the solutions only for the cases~$(a, b) = (0, 1 / 4), (3 / 4, 1 / 4)$, but other combinations also yield similar plots, including a wide range of Pareto solutions.
\Cref{tab:Average computational costs} shows that the new momentum factors are fast enough to compete with the existing ones ($(a, b) = (0, 1/4)$ or $b = a^2/4$) and better than them in some cases.

\begin{table}[htbp]
    \centering
    \ra{1.15}
    \caption{Average computational costs to solve the multi-objective examples}
    \label{tab:Average computational costs}
    \begin{minipage}{.49\hsize}
        \centering
        \subcaption{\cref{eq:JOS1}}
        \begin{tabular}{@{}cccc@{}}
            \toprule
            $a$ & $b$ & Time [\si{\second}] & Iterations \\
            \midrule
        \begin{filecontents*}{data_complexity/JOS1_ab.csv}
0 ,0,6.442,97.0
0 ,1/8,5.158,81.217
0 ,1/4,4.207,65.0
1/6 ,1/144,4.244,67.0
1/6 ,37/288,5.182,82.0
1/6 ,1/4,4.268,66.0
1/4 ,1/64,6.224,99.0
1/4 ,17/128,7.239,113.566
1/4 ,1/4,3.205,51.0
1/2 ,1/16,4.51,72.0
1/2 ,5/32,4.562,71.0
1/2 ,1/4,4.466,70.0
3/4 ,9/64,4.323,67.998
3/4 ,25/128,3.104,49.0
3/4 ,1/4,3.741,47.0
\end{filecontents*}
\csvreader[no head,late after line=\\]{data_complexity/JOS1_ab.csv}
        {1=\a, 2=\b,3=\totaltime,4=\iterationcounts}
            { $\a$ & $\b$ & \totaltime & \iterationcounts}
            \midrule
        \end{tabular}
    \end{minipage}
    \begin{minipage}{.49\hsize}
        \centering
        \subcaption{\cref{eq:JOS1_L1}}
        \begin{tabular}{@{}cccc@{}}
            \toprule
            $a$ & $b$ & Time [\si{\second}] & Iterations \\
            \midrule
        \begin{filecontents*}{data_complexity/JOS1_L1_ab.csv}
0 ,0,10.733,157.512
0 ,1/8,11.054,161.065
0 ,1/4,11.122,161.734
1/6 ,1/144,9.85,141.731
1/6 ,37/288,9.994,144.863
1/6 ,1/4,10.399,150.592
1/4 ,1/64,9.271,135.804
1/4 ,17/128,9.463,137.108
1/4 ,1/4,9.662,139.848
1/2 ,1/16,7.439,109.082
1/2 ,5/32,7.642,110.204
1/2 ,1/4,7.723,111.599
3/4 ,9/64,5.253,77.366
3/4 ,25/128,5.39,79.425
3/4 ,1/4,5.678,82.37
\end{filecontents*}
\csvreader[no head,late after line=\\]{data_complexity/JOS1_L1_ab.csv}
        {1=\a, 2=\b,3=\totaltime,4=\iterationcounts}
            { $\a$ & $\b$ & \totaltime & \iterationcounts}
            \bottomrule
        \end{tabular}
    \end{minipage}
    \begin{minipage}{.49\hsize}
        \centering
        \subcaption{\cref{eq:FDS}}
        \begin{tabular}{@{}cccc@{}}
            \toprule
            $a$ & $b$ & Time [\si{\second}] & Iterations \\
            \midrule
        \begin{filecontents*}{data_complexity/FDS_ab.csv}
0 ,0,29.24,204.438
0 ,1/8,29.797,210.595
0 ,1/4,30.565,214.934
1/6 ,1/144,24.964,174.393
1/6 ,37/288,25.375,177.944
1/6 ,1/4,26.065,182.398
1/4 ,1/64,22.94,159.737
1/4 ,17/128,23.311,162.629
1/4 ,1/4,23.976,166.918
1/2 ,1/16,17.909,122.653
1/2 ,5/32,18.14,123.96
1/2 ,1/4,18.221,125.697
3/4 ,9/64,13.584,94.176
3/4 ,25/128,13.674,94.705
3/4 ,1/4,13.795,94.868
\end{filecontents*}
\csvreader[no head,late after line=\\]{data_complexity/FDS_ab.csv}
        {1=\a, 2=\b,3=\totaltime,4=\iterationcounts}
            { $\a$ & $\b$ & \totaltime & \iterationcounts}
            \bottomrule
        \end{tabular}
    \end{minipage}
    \begin{minipage}{.49\hsize}
        \centering
        \subcaption{\cref{eq:FDS_CONSTRAINED}}
        \begin{tabular}{@{}cccc@{}}
            \toprule
            $a$ & $b$ & Time [\si{\second}] & Iterations \\
            \midrule
        \begin{filecontents*}{data_complexity/FDS_CONSTRAINED_ab.csv}
0 ,0,37.345,259.508
0 ,1/8,37.439,261.522
0 ,1/4,37.94,263.911
1/6 ,1/144,32.463,227.063
1/6 ,37/288,38.265,229.736
1/6 ,1/4,45.661,231.958
1/4 ,1/64,41.434,209.35
1/4 ,17/128,33.664,211.69
1/4 ,1/4,30.772,213.811
1/2 ,1/16,22.92,158.448
1/2 ,5/32,23.1,159.685
1/2 ,1/4,23.539,162.226
3/4 ,9/64,17.092,118.616
3/4 ,25/128,17.123,118.063
3/4 ,1/4,17.115,118.844
\end{filecontents*}
\csvreader[no head,late after line=\\]{data_complexity/FDS_CONSTRAINED_ab.csv}
        {1=\a, 2=\b,3=\totaltime,4=\iterationcounts}
            { $\a$ & $\b$ & \totaltime & \iterationcounts}
            \bottomrule
        \end{tabular}
    \end{minipage}
\end{table}

\subsection{Image deblurring (single-objective)}
\begin{figure}[htpb]
    \centering
    \begin{minipage}[b]{.45\hsize}
        \centering
        \includegraphics[width=0.8\textwidth]{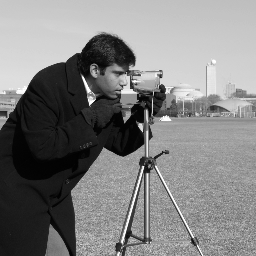}
        \subcaption{Original}
        \label{fig:cameraman:original}
    \end{minipage}
    \begin{minipage}[b]{.45\hsize}
        \centering
        \includegraphics[width=0.8\textwidth]{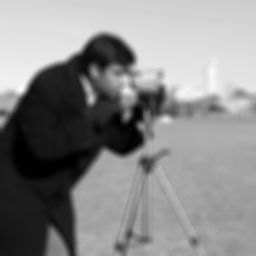}
        \subcaption{Blurred and noisy}
        \label{fig:cameraman:blurred_and_noisy}
    \end{minipage}
    \caption{Deblurring of the cameraman}
    \label{fig:cameraman}
\end{figure}

Since our proposed momentum factor is also new in the single-objective context, we also tackle deblurring the cameraman test image via a single-objective $\ell_2$-$\ell_1$ minimization, inspired by~\cite{Beck2009}.
In detail, as shown in \cref{fig:cameraman}, to a~$256 \times 256$ cameraman test image with each pixel scaled to~$[0,1]$, we generate an observed image by applying a Gaussian blur of size~$9 \times 9$ and standard deviation $4$ and adding a zero-mean white Gaussian noise with standard deviation $10^{-3}$.

Letting~$\theta, B,$ and~$W$ be the observed image, the blur matrix, and the inverse of the Haar wavelet transform, respectively, consider the single-objective problem~\cref{eq:MOP} with~$m = 1$ and
\[ \label{eq:cam_deblur} \tag{CAM-DEBLUR}
    f_1(x) \coloneqq \norm*{BWx - \theta}^2 \eqand g_1(x)=\lambda \norm*{x}_1
,\] 
where~$\lambda \coloneqq 2 \times 10^{-5}$ is a regularization parameter.
Unlike in the previous subsection, we can compute~$\nabla f$'s Lipschitz constant by calculating~$(BW)^\T(BW)$'s eigenvalues using the two-dimensional cosine transform~\citep{Hansen2006}, so we use it constantly as~$\ell$.
Moreover, we use the observed image's Wavelet transform as the initial point.

\begin{figure}[htpb]
    \centering
    \begin{minipage}[b]{.49\hsize}
        \centering
        \includegraphics[width=\linewidth]{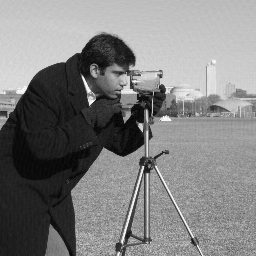}
        \subcaption{$(a, b) = (0 , 1 / 4)$}
    \end{minipage}
    \begin{minipage}[b]{.49\hsize}
        \centering
        \includegraphics[width=\linewidth]{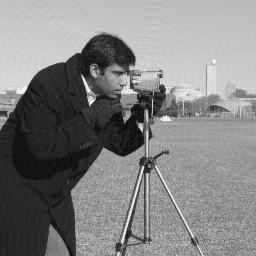}
        \subcaption{$(a, b) = (1 / 2, 1 / 4)$}
    \end{minipage}
    \caption{Deblurred image}
    \label{fig:deblurred}
\end{figure}

\begin{figure}[htpb]
    \centering
    \includegraphics[width=\textwidth]{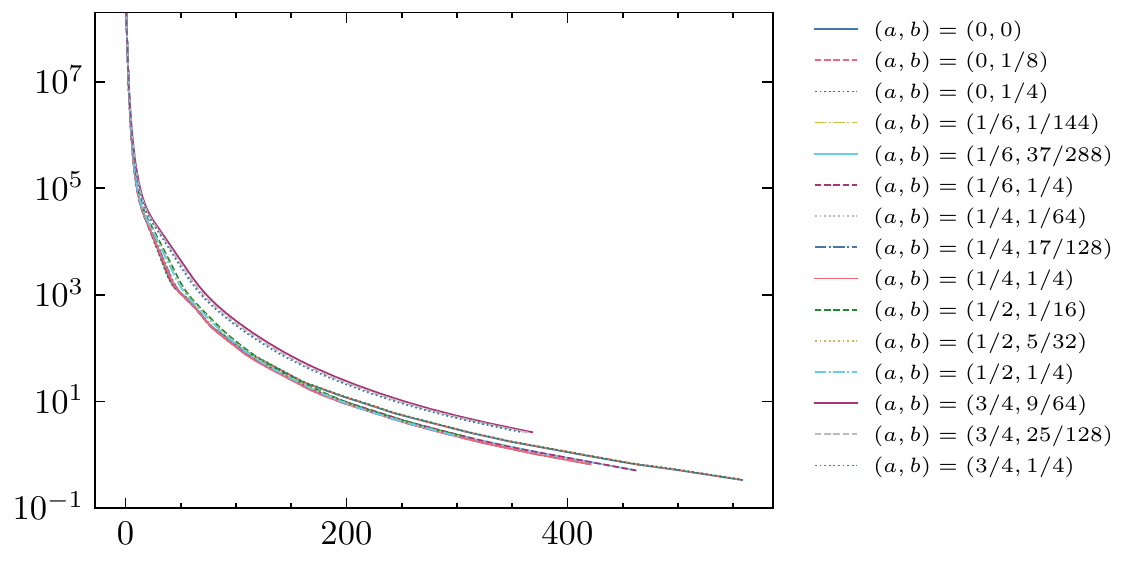}
    \caption{Values of~$u_0(x^k) = F_1(x) - F_1(x^\ast)$, where~$x^\ast$ is the optimal solution estimated from the original image}
    \label{fig:cameraman_plot}
\end{figure}

\Cref{fig:deblurred} shows the reconstructed image from the obtained solution.
Images produced by all hyperparameters are similar, so we present only~$(a, b) = (0, 1 / 4)$ and~$(1 / 2, 1 / 4)$.
Moreover, we summarize the numerical performance in \cref{tab:cam_deblur,fig:cameraman_plot}.
Like the last subsection, this example also suggests that our new momentum factors may occasionally improve the algorithm's performance even for single-objective problems.

\begin{table}[htbp]
    \centering
    \ra{1.15}
    \caption{Computational costs for the image deblurring}
    \label{tab:cam_deblur}
    \begin{tabular}{@{}cccc@{}}
        \toprule
         $a$ & $b$ & Total time [\si{\second}] & Iteration counts \\
        \midrule
    \begin{filecontents*}{data_complexity/cameraman_ab.csv}
0 ,0,85.341,558
0 ,1/8,84.796,558
0 ,1/4,84.75,558
1/6 ,1/144,73.69,460
1/6 ,37/288,73.853,462
1/6 ,1/4,73.903,462
1/4 ,1/64,68.144,421
1/4 ,17/128,67.447,421
1/4 ,1/4,38.452,421
1/2 ,1/16,27.667,305
1/2 ,5/32,27.995,304
1/2 ,1/4,24.649,303
3/4 ,9/64,33.764,368
3/4 ,25/128,35.112,364
3/4 ,1/4,34.05,360
\end{filecontents*}
\csvreader[no head,late after line=\\]{data_complexity/cameraman_ab.csv}
    {1=\a, 2=\b,3=\totaltime,4=\iterationcounts}
        { $\a$ & $\b$ & \totaltime & \iterationcounts}
       \bottomrule
    \end{tabular}
\end{table}

\ifthenelse{\boolean{isMain}}{ 
}{
   {\bibliographystyle{jorsj} 
   \bibliography{library}} 
}

\section{Conclusion} \label{sec: conclusion}
We have generalized the momentum factor of the multi-objective accelerated proximal gradient algorithm~\citep{Tanabe2022a} in a form that is even new in the single-objective context and includes the known FISTA momentum factors~\citep{Beck2009,Chambolle2015}.
Furthermore, with the proposed momentum factor, we proved under reasonable assumptions that the algorithm has an~$O(1/k^2)$ convergence rate and that the iterates converge to Pareto solutions.
Moreover, the numerical results reinforced these theoretical properties and suggested the potential for our new momentum factor to improve the performance.

As we mentioned in \cref{sec: convergence sequence}, our proposed method has the potential to achieve finite-time manifold (active set) identification~\citep{Sun2019} without the assumption of the strong convexity (or its generalizations such as PL conditions or error bounds~\citep{Karimi2016}).
Moreover, we took a single update rule of~$t_k$ for all iterations in this work, but the adaptive change of the strategy in each iteration is conceivable.
It might also be interesting to estimate the Lipschitz constant simultaneously with that change, like in~\cite{Scheinberg2014}.
In addition, an extension to the inexact scheme like~\cite{Villa2013} would be significant.
Those are issues to be addressed in the future.

\ifthenelse{\boolean{isMain}}{ 
}{
   {\bibliographystyle{jorsj} 
   \bibliography{library}} 
}

\section*{Acknowledgements}
This work was supported by the Grant-in-Aid for Scientific Research (C) (21K11769 and 19K11840) and Grant-in-Aid for JSPS Fellows (20J21961) from the Japan Society for the Promotion of Science.

\bibliography{001_library}

\end{document}